\documentclass[10pt,sumlimits,namelimits]{article}

\newcommand{\nc}{\newcommand}
\newcommand{\rnc}{\renewcommand}
\usepackage{amsthm}
\usepackage{amsmath}
\usepackage{amsfonts}
\usepackage{amssymb}
\usepackage{graphicx}
\usepackage{color}
\usepackage[colorlinks=true,urlcolor=blue]{hyperref}
\usepackage[applemac]{inputenc}

\nc{\N}{\mathbb{N}}
\nc{\Z}{\mathbb{Z}}
\nc{\D}{\mathbb{D}}
\nc{\Q}{\mathbb{Q}}
\nc{\R}{\mathbb{R}}
\nc{\C}{\mathbb{C}}
\rnc{\S}{\mathbb{S}}
\nc{\vC}{{\cal C}}
\nc{\vD}{{\cal D}}
\nc{\vF}{{\cal F}}
\nc{\vL}{{\cal L}}
\nc{\vP}{{\cal P}}
\nc{\vS}{{\cal S}}
\nc{\vT}{{\cal T}}
\nc{\vphi}{\varphi}
\nc{\eps}{\varepsilon}
\nc{\dsp}{\displaystyle}
\nc{\ovl}{\overline}
\nc{\udl}{\underline}
\nc{\vlim}{\lim\limits}
\nc{\vlimsup}{\limsup\limits}
\nc{\vliminf}{\liminf\limits}
\nc{\vsup}{\sup\limits}
\nc{\vmax}{\max\limits}
\nc{\vinf}{\inf\limits}
\nc{\vmin}{\min\limits}
\nc{\vint}{\int\limits}
\nc{\vsum}{\sum\limits}
\nc{\vprod}{\prod\limits}
\nc{\vcup}{\bigcup\limits}
\nc{\tends}{\longrightarrow}
\nc{\inj}{\hookrightarrow}
\nc{\weak}{\rightharpoonup}
\nc{\w}{{\textsl w}}
\nc{\wh}{\widehat}
\nc{\loc}{{\rm loc}}
\nc{\rad}{{\rm rad}}
\nc{\supp}{{\rm supp\:}}
\nc{\supess}{{\rm sup\:ess\:}}
\nc{\esssup}{{\rm ess\:sup\:}}
\nc{\Id}{{\rm Id}}
\rnc{\le}{\leqslant}
\rnc{\ge}{\geqslant}
\rnc{\Re}{{\rm Re}}
\rnc{\Im}{{\rm Im}}
\rnc{\det}{{\rm d\acute et\:}}
\nc{\m}{{\rm min}}
\nc{\M}{{\rm max}}
\rnc{\b}{{\rm b}}

\numberwithin{equation}{section}

\newtheorem{thm}{Theorem}[section]
\newtheorem{prop}[thm]{Proposition}
\newtheorem{cor}[thm]{Corollary}
\newtheorem{lem}[thm]{Lemma}

\theoremstyle{definition}
\newtheorem{rmk}[thm]{Remark}
\newtheorem{defi}[thm]{Definition}

\newenvironment{proof*}{\noindent{\bf Proof.}}{\qed}
\newenvironment{vproof}[1]{\noindent{\bf Proof #1}}{\qed}

\title{\Huge \sc Mass Concentration Phenomena for the $\boldmath L^2\unboldmath$-Critical Nonlinear Schrödinger Equation}
\author{\sc Pascal Bégout\footnote{$^{,\dagger}$ Research partially supported by the European network HPRN--CT--2001--00273--HARP ({\it Harmonic analysis and related problems})} and Ana Vargas\footnote{Partially supported by Grant MTM2004--00678 of the MEC (Spain)}}
\date{}

\begin{document}

\maketitle

$$
\begin{array}{cc}
^*\mbox{Laboratoire Jacques-Louis Lions } & \;^\dagger\mbox{Departamento de Matem{\'a}ticas} \\
      \mbox{Université Pierre et Marie Curie } & \mbox{ Facultad de Ciencias} \\
                              \mbox{Boîte Courrier 187 } & \mbox{ Universidad Aut\'onoma de Madrid} \\
                                  \mbox{4, place Jussieu } & \mbox{ Cantoblanco} \\
    \mbox{75252 Paris Cedex 05, FRANCE } & \mbox{ 28049 Madrid, SPAIN}
\bigskip \\
\mbox{
{\footnotesize e-mail\:: }\htmladdnormallink{{\footnotesize\udl{\tt{begout@ann.jussieu.fr}}}}
{mailto:begout@ann.jussieu.fr} } & \mbox{
{\footnotesize e-mail\:: }\htmladdnormallink{{\footnotesize\udl{\tt{ana.vargas@uam.es}}}}
{mailto:ana.vargas@uam.es}}
\bigskip
\end{array}
$$

\begin{abstract}
In this paper, we show that any solution of the nonlinear Schr{\"o}dinger equation $iu_t+\Delta u\pm|u|^\frac{4}{N}u=0,$ which blows up in finite time, satisfies a mass concentration phenomena near the blow-up time. Our proof is essentially based on the Bourgain's one~\cite{MR99f:35184}, which has established this result in the bidimensional spatial case, and on a generalization of Strichartz's inequality, where the bidimensional spatial case was proved by Moyua, Vargas and Vega~\cite{MR1671214}. We also generalize to higher dimensions the results in Keraani~\cite{MR2216444} and  Merle and Vega~\cite{MR1628235}.
\end{abstract}

{\let\thefootnote\relax}
{\let\thefootnote\relax}

\tableofcontents

\baselineskip .7cm

\section{Introduction and main results}
\label{introduction}

Let $\gamma\in\R\setminus\{0\}$ and let $0\le\alpha\le\frac{4}{N}.$ It is well-known that for any $u_0\in L^2(\R^N),$ there exists a unique maximal solution
$$
u\in C((-T_\m,T_\M);L^2(\R^N))\cap L^\frac{4(\alpha+2)}{N\alpha}_\loc((-T_\m,T_\M);L^{\alpha+2}(\R^N)),
$$
of
\begin{gather}
 \left\{
  \begin{split}
   \label{nls*}
    i\frac{\partial u}{\partial t}+\Delta u+\gamma|u|^\alpha u= 0, & \; (t,x)\in(-T_\m,T_\M)\times\R^N, \\
                                                                                          u(0)= u_0, & \mbox{ in }\R^N, \\
  \end{split}
 \right.
\end{gather}
satisfying the conservation of charge, that is for any $t\in(-T_\m,T_\M)$, $\|u(t)\|_{L^2(\R^N)}=\|u_0\|_{L^2(\R^N)}.$
The solution $u$ also satisfies the following Duhamel's formula
\begin{gather}
\label{inteq}
\forall t\in(-T_\m,T_\M),\;
u(t)=\vT(t)u_0+i\gamma\vint_0^{\;t}(\vT(t-s)\{|u|^\alpha u\})(s)ds,
\end{gather}
where we design by $(\vT(t))_{t \in \R}$ the group of isometries $(e^{it\Delta})_{t\in\R}$ generated by $i\Delta$ on
$L^2(\R^N;\C).$ Moreover $u$ is maximal in the following sense. If $\alpha<\frac{4}{N}$ then $T_\M=T_\m=\infty,$ if
$\alpha=\frac{4}{N}$ and if $T_\M<\infty$ then
$$
\|u\|_{L^\frac{2(N+2)}{N}((0,T_\M);L^\frac{2(N+2)}{N}(\R^N))}=\infty,
$$
and if $\alpha=\frac{4}{N}$ and $T_\m<\infty$ then $\|u\|_{L^\frac{2(N+2)}{N}((-T_\m,0);L^\frac{2(N+2)}{N}(\R^N))}=\infty$ (see Cazenave and Weissler~\cite{MR91a:35149} and Tsutsumi~\cite{MR915266}, also
Cazenave~\cite{MR2002047}, Corollary~4.6.5 and Section~4.7). Now, assume that $\alpha=\frac{4}{N}.$ It is well-known that if $\|u_0\|_{L^2}$ is small enough then $T_\M=T_\m=\infty,$ whereas if $\gamma>0$ then there exists some $u_0\in L^2(\R^N)$ such that $T_\M<\infty$ and $T_\m<\infty.$ For example, it is sufficient to
choose $u_0=\lambda\vphi,$ where $\vphi\in H^1(\R^N)\cap L^2(|x|^2;dx),$ $\vphi\not\equiv0,$ and where $\lambda>0$ is large enough (Glassey~\cite{MR57:842}, Vlasov, Petrischev and Talanov~\cite{vpt}, Cazenave and Weissler~\cite{MR91a:35149}).
\bigskip
\\
In the case $\gamma>0,$ when blow-up in finite time occurs, a mass concentration phenomena was observed near the blow-up time (see Theorem~2 in Merle and Tsutsumi~\cite{MR91e:35194} and Theorem~6.6.7 in Cazenave~\cite{MR2002047}), under the conditions that $u_0\in H^1(\R^N)$ is spherically symmetric, $N\ge2$ and $\gamma>0.$ Theorem~6.6.7 in Cazenave~\cite{MR2002047} asserts that if $T_\M<\infty$ for a solution $u$ of equation (\ref{nls}) below, then for any
$\eps\in\left(0,\frac{1}{2}\right),$
\begin{gather}
\label{intro0}
\vliminf_{t\nearrow T_\M}\int_{B(0,(T_\M-t)^{\frac{1}{2}-\eps})}|u(t,x)|^2dx\ge \|Q\|_{L^2(\R^N)}^2,
\end{gather}
where $Q$ is the ground state, {\it i.e.} the unique positive solution of $-\Delta Q+Q=|Q|^\frac{4}{N}Q$ (see Merle and Tsutsumi~\cite{MR91e:35194}, Tsutsumi~\cite{MR915266}). The proof uses the conservation of energy and the compactness property of radially symmetric functions lying in $H^1(\R^N).$ The spherical symmetry assumption was relaxed by Nawa~\cite{MR1172417}; see also Hmidi and Keraani~\cite{MR2180464}. Later, it was proved that for data in $H^s,$ for some $s<1,$ \eqref{intro0} holds. This was proved by  Colliander, Raynor, Sulem and Wright~\cite{MR2150890} for dimension $2,$ and extended by Tzirakis~\cite{MR2213400} to dimension $1$ and by Visan and Zhang~\cite{MR2318374} to general dimension.
\bigskip
\\
In Bourgain~\cite{MR99f:35184}, a mass concentration phenomena, estimate (\ref{intro1}) below, is obtained for any $u_0\in L^2(\R^2),$ $\gamma\neq0,$ but in spatial dimension $N=2.$ Consider solutions of the following critical nonlinear Schr{\"o}dinger equation,
\begin{gather}
 \left\{
  \begin{split}
   \label{nls}
    i\frac{\partial u}{\partial t}+\Delta u+\gamma|u|^\frac{4}{N}u= 0, & \; (t,x)\in(-T_\m,T_\M)\times\R^N, \\
                                                                                                u(0)= u_0, & \mbox{ in }\R^N, \\
  \end{split}
 \right.
\end{gather}
where $\gamma\in\R\setminus\{0\}$ is a given parameter. Bourgain showed, in the case $N=2$ (see Theorem~1 in~\cite{MR99f:35184}), that if $u\in C((-T_\m,T_\M);L^2(\R^2))$ is a solution of (\ref{nls}) with initial data $u_0\in L^2(\R^2)$ which blows-up in finite time $T_\M<\infty,$ then
\begin{gather}
\label{intro1}
\vlimsup_{t\nearrow T_\M}\sup_{c\in\R^N}\int_{B(c,C(T_\M-t)^\frac{1}{2})}|u(t,x)|^2dx\ge\eps,
\end{gather}
where the constants $C$ and $\eps$ depend continuously and only on $\|u_0\|_{L^2}$ and $|\gamma|.$ The proof is based on a refinement of Strichartz's inequality for $N=2,$ due to Moyua, Vargas and Vega (see Theorem~4.2 and Lemma~4.4 in~\cite{MR1671214}).
\bigskip
\\
Very recently, Keraani~\cite{MR2216444} showed for $N\in\{1,2\}$ that there is some $\delta_0>0,$ such that, under the same assumptions, if in addition $\|u_0\|_{L^2}<\sqrt2\delta_0$ then for any $\lambda(t)>0$ such that $\lambda(t)\xrightarrow{t\nearrow T_\M}\infty,$
\begin{gather}
\label{kera}
\vliminf_{t\nearrow T_\M}\sup_{c\in\R^N}
\int_{B(c,\lambda(t)(T_\M-t)^\frac{1}{2})}|u(t,x)|^2dx\ge\delta_0^2.
\end{gather}
Keraani's proof uses a linear profile decomposition that was shown in dimension $N=2$ by Merle and Vega~\cite{MR1628235} and in dimension $N=1$ by Carles and Keraani~\cite{MR2247881} (see Theorem~\ref{thmlinearprofiles} below for the precise statement). The proofs of the decompositions are based on the above mentioned
refinement of Strichartz's inequality by Moyua, Vargas and Vega and another one for the case $N=1$ observed by Carles and Keraani~\cite{MR2247881}. In this paper, we generalize the refinement of Strichartz's inequality (see Theorem~\ref{stric3} below) in order to establish the higher dimensional versions of all these results. Our proofs (namely, those of Theorem~\ref{stric1} and Lemma~\ref{lemg}) rely on the restriction theorems for paraboloids
proved by Tao~\cite{MR2033842}. There is another minor technical point, because the Strichartz's exponent $\frac{2N+4}N,$ is not a natural number when the dimension $N\ge 3,$ except $N=4.$ We have to deal with this little inconvenience which did not appeared in $N\in\{1,2\}.$
\bigskip
\\
This paper is organized as follows. At the end of this section, we state the main results (Theorems~\ref{thmmass} and \ref{stric3}) and give some notations which will be used throughout this paper. Section~\ref{refinement} is devoted to the proof of the refinement of Strichartz's inequality (Theorems~\ref{stric1}--\ref{stric3}). In Section~\ref{preli}, we establish some preliminary results in order to prove a mass concentration result in Section~\ref{mass} (Proposition~\ref{propmass}). We prove Theorem~\ref{thmmass} in Section~\ref{mass}. Finally, Section~\ref{further} is devoted to the generalization to higher dimensions of the results by Keraani~\cite{MR2216444} and Merle and Vega~\cite{MR1628235}.
\bigskip
\\
Throughout this paper, we use the following notation. For $1\le p\le\infty,$ $p'$ denotes the conjugate of $p$ defined by $\frac{1}{p}+\frac{1}{p'}=1;$ $L^p(\R^N)=L^p(\R^N;\C)$ is the usual Lebesgue space. The Laplacian in $\R^N$ is written $\Delta=\sum\limits_{j=1}^N\frac{\partial^2}{\partial x^2_j}$ and $\frac{\partial u}{\partial t}=u_t$ is the time derivative of the complex-valued function $u.$ For $c\in\R^N$ and $R\in(0,\infty),$ we denote by $B(c,R)=\{x\in \R^N;\; |x-c|<R\}$ the open ball of $\R^N$ of center $c$ and radius $R.$ We design by $\vC$ the set of half--closed cubes in $\R^N$. So $\tau\in\vC$ if and only if there exist $(a_1,\ldots,a_N)\in\R^N$ and $R>0$ such that $\tau=\vprod_{j=1}^N[a_j,a_j+R).$ The length of a side of $\tau\in\vC$ is written $\ell(\tau)=R.$ Given $A\subset\R^N,$ we denote by $|A|$ its Lebesgue measure. Let $j,k\in\N$ with $j<k.$ Then we denote $[\![j,k]\!]=[j,k]\cap\N.$ We denote by $\vF$ the Fourier transform in $\R^N$ defined by
\footnote{with this definition of the Fourier transform, $\|\vF u\|_{L^2}=\|\vF^{-1}u\|_{L^2}=\|u\|_{L^2},$
                 $\vF^{-1}\vF=\vF\vF^{-1}=\Id_{L^2},$ $\vF(u*v)=\vF u\vF v$ and $\vF^{-1}(u*v)=\vF^{-1}u\vF^{-1}v.$}
$\wh u(\xi)=\vF u(\xi)=\dsp\int_{\R^N}e^{-2i\pi x.\xi}u(x)dx,$ and by $\vF^{-1}$ its inverse given by
$\vF^{-1}u(x)=\dsp\int_{\R^N}e^{2i\pi\xi.x}u(\xi)d\xi.$ $C$ are auxiliary positive constants and $C(a_1, a_2,\dots,a_n)$ indicates that the constant $C$ depends only on positive parameters $a_1,a_2,\dots, a_n$ and that the dependence is continuous.
\bigskip
\\
Finally, we recall the Strichartz's estimates (Stein--Tomas Theorem) (see Stein~\cite{MR864375}, Strichartz~\cite{MR0512086} and
Tomas~\cite{MR0358216}). Let $I\subseteq\R$ be an interval, let $t_0\in\ovl I$ and let $\gamma\in\C.$ Set for any $t\in I,$
$\Phi_u(t)=i\gamma\dsp\int_{t_0}^t(\vT(t-s)\{|u|^\frac{4}{N}u\})(s)ds.$ Then we have
\begin{gather}
\label{stri1}
\|\vT(\:.\:)u_0\|_{L^\frac{2(N+2)}{N}(\R\times\R^N)}\le C_0\|u_0\|_{L^2(\R^N)}, \medskip \\
\label{stri2}
\|\Phi_u\|_{L^\frac{2(N+2)}{N}(I\times\R^N)}\le C_1\|u\|_{L^\frac{2(N+2)}{N}(I\times\R^N)}^\frac{N+4}{N},
\end{gather}
where $C_0=C_0(N)>0$ and $C_1=C_1(N,|\gamma|)>0.$ For more details, see Ginibre and Velo~\cite{MR87b:35150} (Lemma~3.1) and Cazenave and Weissler~\cite{MR91a:35149} (Lemma~3.1), also Cazenave~\cite{MR2002047} (Theorem~2.3.3). The main results of this paper are the following.

\begin{thm}
\label{thmmass}
Let $\gamma\in\R\setminus\{0\},$ let $u_0\in L^2(\R^N)\setminus\{0\}$ and let
$$
u\in C((-T_\m,T_\M);L^2(\R^N))\cap L_\loc^\frac{2(N+2)}{N}((-T_\m,T_\M);L^\frac{2(N+2)}{N}(\R^N))
$$
be the maximal solution of $(\ref{nls})$ such that $u(0)=u_0.$ There exists $\eps=\eps(\|u_0\|_{L^2},N,|\gamma|)>0$  satisfying the following property. If $T_\M<\infty$ then
\begin{gather*}
\vlimsup_{t\nearrow T_\M}\sup_{c\in\R^N}\int_{B(c,(T_\M-t)^\frac{1}{2})}|u(t,x)|^2dx\ge\eps,
\end{gather*}
and if $T_\m<\infty$ then
\begin{gather*}
\vlimsup_{t\searrow-T_\m}\sup_{c\in\R^N}\int_{B(c,(T_\m+t)^\frac{1}{2})}|u(t,x)|^2dx\ge\eps.
\end{gather*}
\end{thm}
\noindent
By keeping track of the constants through the proofs, it can be shown that $\eps=C(N,|\gamma|)\|u_0\|_{L^2}^{-m}$ for some $m>0$ (this was pointed out by Colliander). Notice that no hypothesis on the attractivity on the nonlinearity (that is on the $\gamma$'s sign), on the spatial dimension $N$ and on the smoothness on the initial data $u_0$ are made.
\bigskip
\\
For each $j\in\Z,$ we break up $\R^N$ into dyadic cubes $\tau_k^j=\vprod_{m=1}^N[k_m 2^{-j},(k_m+1) 2^{-j}),$ where $k=(k_1,\dots,k_N)\in\Z^N$ with $\ell(\tau_k^j)=2^{-j}.$ Define $f_k^j(x)=f\chi_{\tau_k^j}(x).$ Let $1\le p<\infty$ and let $1\le q<\infty.$ We define the spaces
$$
X_{p,q}=\left\{f\in L^p_\loc(\R^N);\; \|f\|_{X_{p,q}}<\infty\right\},
$$
where
$$
\|f\|_{X_{p,q}}
=\left[\sum_{j\in\Z}2^{j\frac{N}{2}\frac{2-p}p q}\sum_{k\in\Z^N}\|f_k^j\|_{L^p(\R^N)}^q\right]^\frac{1}{q}.
$$
Then $(X_{p,q},\|\;.\;\|_{X_{p,q}})$ is a Banach space and the set of functions  $f\in L^\infty(\R^N)$  with compact support is dense in $X_{p,q}$ for the norm $\|\;.\;\|_{X_{p,q}}.$
\bigskip
\\
We prove the following improvement of Strichartz's (Stein--Tomas's) inequality.

\begin{thm}
\label{stric1}
Let $q=\frac{2(N+2)}{N}$ and $1<p<2$ be such that $\frac{1}{p'}>\frac{N+3}{N+1}\frac{1}{q}.$ For every function $g$ such that $g\in X_{p,q}$ or $\wh g\in X_{p,q},$ we have
\begin{gather}
\label{stric1-1}
\|\vT(\:.\:)g\|_{L^q(\R^{N+1})}\le C\min\left\{\|g\|_{X_{p,q}},\|\wh g\|_{X_{p,q}}\right\},
\end{gather}
where $C=C(N,p).$
\end{thm}

\begin{thm}
\label{stric2}
Let $q>2$ and let $1<p<2.$  Then there exists $\mu\in\left(0,\frac{1}{p}\right)$ such that for every function
$f\in L^2(\R^N),$ we have
\begin{gather}
\|f\|_{X_{p,q}}\le C\left[\sup_{(j,k)\in\Z\times\Z^N}2^{j\frac N2(2-p)}
\int_{\tau_k^j}|f(x)|^pdx\right]^\mu\|f\|_{L^2(\R^N)}^{1-\mu p}\le C\|f\|_{L^2(\R^N)},
\end{gather}
where $C=C(p,q)$ and $\mu=\mu(p,q).$ In particular, $L^2(\R^N)\inj
X_{p,q}.$ Moreover, $L^2(\R^N)\neq X_{p,q}.$
\end{thm}
\noindent
As a corollary we obtain the following improvement of Strichartz's (Stein--Tomas's) inequality.

\begin{thm}
\label{stric3}
Let $q=\frac{2(N+2)}{N}$ and let $p<2$ be such that $\frac{1}{p'}>\frac{N+3}{N+1}\frac{1}{q}.$ Then, there exists
$\mu\in\left(0,\frac{1}{p}\right)$ such that for every function $g\in L^2(\R^N),$ we have
\begin{gather}
\label{stric2-1}
\|\vT(\:.\:)g\|_{L^q(\R^{N+1})}\le C\left[\sup_{(j,k)\in\Z\times\Z^N}2^{j\frac N2(2-p)}
\int_{\tau_k^j}|\wh g(\xi)|^pd\xi\right]^{\mu}\|g\|_{L^2(\R^N)}^{1-\mu p}\le C\|g\|_{L^2(\R^N)},
\end{gather}
where $C=C(N,p)$ and $\mu=\mu(N,p).$
\end{thm}

\begin{rmk}[See Bourgain~\cite{MR99f:35184}, p.262--263]
\label{rmkfurther}
By H\"older's inequality, if $1<p<2$ then for any $(j,k)\in\Z\times\Z^N,$
$$
\left[2^{j\frac N2(2-p)}\int_{\tau_k^j}|\wh g(\xi)|^pd\xi \right]^{1/p}\le \left[2^{j\frac N2}
\int_{\tau_k^j}|\wh g(\xi)|d\xi \right]^\theta\| \wh g\|_{L^2(\R^N)}^{1-\theta}\le \|g\|_{B_{2,\infty}^0}^\theta
\| \wh g\|_{L^2(\R^N)}^{1-\theta},
$$
for some $0<\theta<1.$ Therefore, it follows from our Strichartz's refinement, Theorem~\ref{stric3}, that the following holds.
$$
\forall M>0,\; \exists\eta>0\mbox{ such that if } \|u_0\|_{L^2}\le M \mbox{ and } \|u_0\|_{B^0_{2,\infty}}<\eta
\mbox{ then } T_\M=T_\m=\infty,
$$
where $u$ is the corresponding solution of ~\eqref{nls}. Furthermore, $u\in
L^\frac{2(N+2)}{N}(\R;L^\frac{2(N+2)}{N}(\R^N))$ and there exists a scattering state in $L^2(\R^N).$ The same result holds if the condition $\|u_0\|_{B^0_{2,\infty}}<\eta$ is replaced by
$$
\sup_{(j,k)\in\Z\times\Z^N}2^{j\frac N2(2-p)}\int_{\tau_k^j}|u_0(x)|^pdx<\eta',
$$
for a suitable $\eta^\prime.$
\end{rmk}

Very recently, Rogers and Vargas~\cite{MR2264250} have proved, for the non--elliptic cubic Schr\"odinger equation $i\partial_tu+\partial^2_{x_1}u-\partial^2_{x_2}u+\gamma |u|^2u=0$ in dimension 2, some results analogous to Theorems~\ref{thmmass}, \ref{stric1}, \ref{stric2} and \ref{stric3}.

\section{Strichartz's refinement}
\label{refinement}

We recall that $\vT(t)g=K_t*g,$ where $K_t(x)=(4\pi it)^{-\frac{N}{2}}e^{i\frac{|x|^2}{4t}}$ and that
$\wh{K_t}(\xi)=e^{-i4\pi^2|\xi|^2t}.$ Using that for any $g\in L^2(\R^N),$ $\vT(t)g=\vF^{-1}(\wh{K_t}\wh g)$ we have,
\begin{equation}
\label{demlemg1}
(\vT(t)g)(x)=\vint_{\R^N}e^{2i\pi\left(x.\xi-2\pi t|\xi|^2\right)}\wh g(\xi)d\xi.
\end{equation}
Let $S=\left\{(\tau,\xi)\in\R\times\R^N;\; \tau=-2\pi|\xi|^2\right\},$ let $d\sigma(|\xi|^2,\xi)=d\xi$ and
let $f$ be defined on $S$ by
$f(\tau,\xi)=f(-2\pi|\xi|^2,\xi)=\wh{g}(\xi).$ Then,
\begin{gather}
 \begin{split}
   \label{adjrest}
      & (\vT(t)g)(x)=\vint_{\R^N}f(-2\pi|\xi|^2,\xi)e^{2i\pi(x.\xi-2\pi t|\xi|^2)}d\xi \medskip \\
   = & \iint\limits_{\!S}f(\tau,\xi)e^{2i\pi(t\tau+x.\xi)}d\sigma(\tau,\xi)=\vF^{-1}(fd\sigma)(t,x).
 \end{split}
\end{gather}
Our main tool will be the following bilinear restriction estimate proved by Tao~\cite{MR2033842}. We adapt the statements to our notation using the equivalence \eqref{adjrest}.

\begin{thm}[Theorem~1.1 in~\cite{MR2033842}]
Let $Q,$ $Q'$ be cubes of sidelength $1$ in $\R^N$ such that
$$
\min\{d(x,y);\; x\in Q,\; y\in Q'\}\sim1
$$
and let $\wh{f},$ $\wh{g}$ functions respectively supported in $Q$ and $Q^\prime.$ Then for any $r>\frac{N+3}{N+1}$ and $p\ge2,$ we have
$$
\|\vT(\:.\:)f \vT(\:.\:)g\|_{L^r(\R^{N+1})}\le C\|\wh{f}\|_{L^p(Q)}\|\wh{g}\|_{L^p(Q')},
$$
with a constant $C$ independent of $f,$ $g,$ $Q$ and $Q'.$
\end{thm}

By interpolation with the trivial estimate
$$
\|\vT(\:.\:)f \vT(\:.\:)g\|_{L^\infty(\R^{N+1})}\le C\|\wh{f}\|_{L^1(Q)}\|\wh{g}\|_{L^1(Q')}
\le C\|\wh{f}\|_{L^p(Q)}\|\wh{g}\|_{L^p(Q')},
$$
for any $p\ge1,$ one obtains the following result.

\begin{thm}[\cite{MR2033842}]
Let $Q,$ $Q'$ be cubes of sidelength $1$ in $\R^N$ such that
$$
\min\{d(x,y);\; x\in Q,\; y\in Q'\}\sim1
$$
and $\wh{f},$ $\wh{g}$ functions respectively supported in $Q$ and $Q^\prime.$ Then for any $r>\frac{N+3}{N+1}$ and for all $p$ such that $\frac2{p'}>\frac{N+3}{N+1}\frac1r,$ we have
$$
\|\vT(\:.\:)f \vT(\:.\:)g\|_{L^r(\R^{N+1})}\le C\|\wh{f}\|_{L^p(\R^N)}\|\wh{g}\|_{L^p(\R^N)},
$$
with a constant $C$ independent of $f,$ $g,$ $Q$ and $Q'.$
\end{thm}
By rescaling and taking $r=\frac{N+2}N,$ we obtain the following.

\begin{cor}
\label{tao}
Let $\tau,$ $\tau'$ be cubes of sidelength $2^{-j}$ such that
$$
\min\{d(x,y);\; x\in \tau,\; y\in \tau'\}\sim2^{-j}
$$
and $\wh{f},$ $\wh{g}$ functions respectively supported in $\tau$ and $\tau'.$ Then for $r=\frac{N+2}N$ and for any $p$ such that $\frac2{p'}>\frac{N+3}{N+1}\frac1r,$ we have
$$
\|\vT(\:.\:)f \vT(\:.\:)g\|_{L^r(\R^{N+1})}\le C2^{jN\frac{2-p}{p}}\|\wh{f}\|_{L^p(\R^N)}\|\wh{g}\|_{L^p(\R^N)},
$$
with a constant $C$ independent of $f,$ $g,$ $\tau$ and $\tau'.$
\end{cor}

We will need to use the orthogonality of functions with disjoint support. More precisely, the following lemma, a proof of which can be found, for instance, in Tao, Vargas, Vega~\cite{MR1625056}, Lemma~6.1.

\begin{lem}
\label{lem}
Let $(R_k)_{k\in\Z}$ be a collection of rectangles in frequency space and $c>0,$ such that the dilates $(1+c)R_k$ are almost disjoint $($i.e. $\sum_k \chi_{(1+c)R_{k}}\le C),$ and suppose that $(f_k)_{k\in\Z}$ is a collection of functions whose Fourier transforms are supported on $R_k.$ Then for all $1\le p\le\infty,$
we have
$$
\|\sum_{k\in\Z}f_k\|_{L^p(\R^N)}\le C(N,c)\left(\:\sum_{k\in\Z} \|f_k\|_{L^p(\R^N)}^{p^*}\right)^\frac{1}{p^*},
$$
where $p^*=\min(p, p^\prime).$
\end{lem}

\begin{vproof}{of Theorem \ref{stric1}.}
We set $r=\frac q2=\frac{N+2}N.$ We first consider the case where $\wh g\in X_{p,q}.$ We can assume that the support of $\wh{g}$ is contained in the unit square. The general result follows by scaling and density. For each $j\in\Z,$ we decompose $\R^N$ into dyadic cubes $\tau_k^j$ of sidelength $2^{-j}.$ Given a dyadic cube $\tau_k^j$ we will say that it is the ``parent'' of the $2^N$ dyadic cubes of sidelength $2^{-j-1}$ contained in it. We write $\tau_k^j\sim\tau_{k'}^j$ if $\tau_k^j$, $\tau_{k'}^j$ are not adjacent but have adjacent parents. For each $j\ge0,$ write $g=\sum g_k^j$ where $\wh{g}_k^j(\xi)=\wh{g}\chi_{\tau_k^j}(\xi).$ Denote by $\Gamma$ the diagonal of $\R^N\times\R^N,$
$\Gamma=\{(x,x);\;x\in\R^N\}.$ We have the following decomposition (of Whitney type) of $\R^N\times\R^N\setminus\Gamma$ (see Figure~\ref{1}),
$$
(\R^N\times\R^N)\setminus\Gamma=\bigcup_j\bigcup_{k,k';\;\tau_k^j\sim\tau_{k'}^j}\tau_k^j\times\tau_{k'}^j.
$$
\begin{figure}[htbp]
\centering
\includegraphics[ext=jpg,hiresbb=true,width=3in]{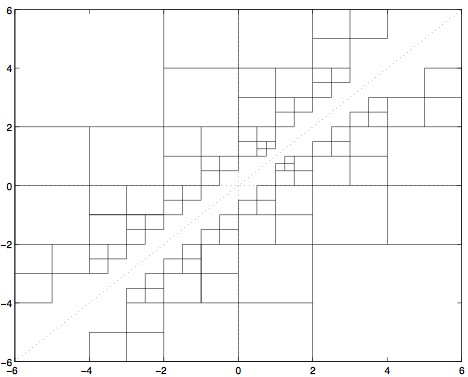}
\caption{$\R^N\times \R^N$}
\label{1}
\end{figure}

Thus,
\begin{gather*}
 \begin{split}
     & \vT(t)g(x)\,\vT(t)g(x)=\vint_{\R^N}\vint_{\R^N}e^{2i\pi\left(x.\xi-2\pi t|\xi|^2\right)}\wh g(\xi)
         e^{2i\pi\left(x.\eta-2\pi t|\eta|^2\right)}\wh g(\eta)d\xi d\eta \medskip \\
  = & \sum_{j}\sum_{k}\sum_{k';\tau^j_k\sim\tau^j_{k'}}\vint\vint_{\tau_k^j\times\tau_{k'}^j}
         e^{2i\pi\left(x.\xi-2\pi t|\xi|^2\right)}\wh g(\xi)e^{2i\pi\left(x.\eta-2\pi t|\eta|^2\right)}\wh g(\eta)d\xi d\eta
         \medskip \\
  = & \sum_{j}\sum_{k}\sum_{k';\tau^j_k\sim\tau^j_{k'}}\vT(t)g_k^j\,\vT(t)g_{k'}^j
 \end{split}
\end{gather*}
(see also Tao, Vargas and Vega~\cite{MR1625056}). Thus,
$$
\|\vT(\:.\:)g\|_{L^{2r}(\R^{N+1})}^{2}=\|\vT(\:.\:)g \vT(\:.\:)g\|_{L^r(\R^{N+1})}
=\| \sum_j \sum_{\underset
{\tau^j_k\sim\tau^j_{k'}}{k,k':}}\vT(\:.\:)g_k^j\vT(\:.\:)
g_{k'}^j\|_{L^r(\R^{N+1})}.
$$
For each $k=(k_1,k_2,\dots, k_N),$ the support of the $(N+1)$-dimensional Fourier transform of $\vT(\:.\:)g_k^j$ is
contained in the set $\tilde\tau_k^j=\{(-2\pi|\xi|^2,\xi);\;\xi\in\tau_k^j\}.$ Hence the support of the Fourier transform of $\vT(\:.\:)g_k^j \vT(\:.\:) g_{k'}^j$ is contained in $\tilde\tau_k^j+\tilde\tau_{k'}^j=\{(-2\pi(|\xi|^2+|\xi'|^2),\xi+\xi');\;
\xi\in\tau_k^j,\;\xi'\in\tau_{k'}^j\}.$ Using the identity $|\xi|^2+|\xi'|^2= \frac12|\xi+\xi'|^2+\frac12|\xi-\xi'|^2$ we see that
$\tilde\tau_k^j+\tilde\tau_{k'}^j$ is contained in the set $H_{j,k}=\{(a,b)\in\R^N\times\R:\; |a-2^{-j+1}k|\le C2^{-j}, \;
2^{-2j}\le-|a|^2-\frac b\pi\le 3 N2^{-2j}\}.$ Note that,
$$
\sum_j \sum_k\sum_{k';\;\tau^j_k\sim\tau^j_{k'}}\chi_{H_{j,k}}\le C(N).
$$
Hence, the functions $\vT(\:.\:)g_k^j\vT(\:.\:)g_{k'}^j$ are almost orthogonal in $L^2(\R^{N+1}).$ A similar orthogonality condition was the key in the proof of the $L^4$--boundedness of the Bochner--Riesz multipliers given by C\'{o}rdoba~\cite{MR0447949}, see also Tao, Vargas and Vega~\cite{MR1625056}, and implicitly appears in Bourgain~\cite{MR1122623}, Moyua, Vargas and Vega~\cite{MR1413873,MR1671214}. But we need something more, since we are not working in $L^2$ and we want to apply Lemma~\ref{lem}. For $M=2 [\ln(N+1)],$ we decompose each $\tau_j^k$ into dyadic subcubes of sidelength $2^{-j-M}.$ Consequently, we have a corresponding decomposition of $\tau_k^j\times\tau_{k'}^j$ and of $\R^N\times\R^N,$ as follows~: set $\mathcal D$ the family of multi-indices $(m,m',\ell)\in\Z^N\times\Z^N\times\Z,$ so that, there exists some $\tau_k^{\ell-M}$ and $\tau_{k'}^{\ell-M}$ with $\tau_{m}^\ell\subset\tau_{k}^{\ell-M},$ $\tau_{m'}^\ell\subset\tau_{k'}^{\ell-M}$ and $\tau_k^{\ell-M}\sim\tau_{k'}^{\ell-M}$
$(j=\ell-M).$ Then,
$$
(\R^N\times\R^N)\setminus\Gamma=\bigcup_\vD\tau_m^\ell\times\tau_{m'}^\ell.
$$
Hence,
$$
\|\vT(\:.\:)g\|_{L^{2r}(\R^{N+1})}^2=\|\vT(\:.\:)g \vT(\:.\:)g\|_{L^r(\R^{N+1})}
=\| \sum_\vD\vT(\:.\:)g_m^\ell\vT(\:.\:)
g_{m'}^\ell\|_{L^r(\R^{N+1})}.
$$
Notice that if $(m,m',\ell)\in\mathcal D,$ then the distance between $\tau_m^\ell$ and $\tau_{m'}^\ell$ is bigger than
$2^{-\ell+M}\ge N 2^{-\ell},$ and smaller than $\sqrt N2^{-\ell+M}.$ We {\bf claim} that there are rectangles
$R_{m,m',\ell},$ and $c=c(N),$ so that $\tilde\tau_m^\ell\times\tilde\tau_{m'}^\ell\subset R_{m,m',\ell}$ and $\sum_\vD\chi_{(1+c)R_{m,m',\ell}}\le C(N).$ We postpone the proof of this claim to the end of the proof. Assuming that it holds, and by Lemma~\ref{lem}, since $r<2,$ we have
$$
\|\sum_\vD \vT(\:.\:)g_m^\ell\vT(\:.\:) g_{m'}^\ell\|_{L^r(\R^{N+1})}\le C(N)
\left[\sum_\vD\|\vT(\:.\:)\,g_m^\ell\vT(\:.\:)g_{m'}^\ell\|_{L^r(\R^{N+1})}^r\right]^\frac{1}{r}.
$$
Now use Corollary~\ref{tao} to estimate
\begin{align*}
     & \left[\sum_\vD\|\vT(\:.\:)\,g_m^\ell\vT(\:.\:)g_{m'}^\ell\|_{L^r(\R^{N+1})}^r\right]^\frac{1}{r} \medskip \\
\le & \;C(N,p)\left[\sum_\ell\sum_{m}\sum_{m';(m,m',\ell)\in\vD}2^{\ell Nr\frac{2-p}{p}}\|\wh{g}_m^\ell\|_{L^p(\R^N)}^r
           \|\wh{g}_{m'}^\ell\|_{L^p(\R^N)}^r\right]^\frac{1}{r}.
\end{align*}
Now, for each $(m,\ell)$ there are at most $4^N2^{MN}$ indices
$m'$ such that $(m,m',\ell)\in\vD.$ Hence,
\begin{gather*}
\left[\sum_\ell\sum_{m}\sum_{m';(m,m',\ell)\in\vD}  2^{\ell Nr\frac{2-p}{p}}\|\wh{g}_m^\ell\|_{L^p(\R^N)}^r
\|\wh{g}_{m'}^\ell\|_{L^p(\R^N)}^r\right]^\frac{1}{r}
\le C(N)\left[\sum_\ell\sum_m  2^{\ell Nr\frac{2-p}{p}}
\|\wh{g}_m^\ell\|_{L^p(\R^N)}^{2r}\right]^\frac{1}{r}.
\end{gather*}

\begin{figure}[htbp]
\centering
\includegraphics[ext=jpg,hiresbb=true,width=3in]{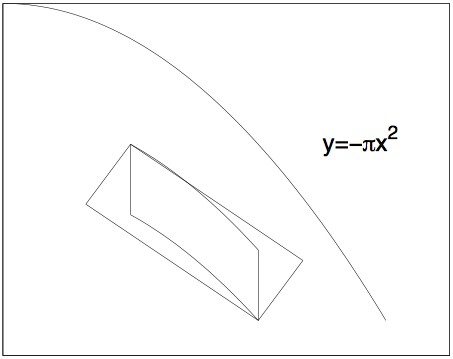}
\caption{$H_{m,m',\ell}\subset R_{m,m',\ell}$ }
\label{rectangle}
\end{figure}
We still have to justify the claim.  Assume, for the sake of simplicity that
$$
\tau_m^\ell\times\tau_{m'}^\ell\subset\{(x_1,x_2,\ldots,x_N)\in\R^N;\; \forall i\in[\![1,N]\!],\; x_i\ge0\}.
$$
Then $\tilde\tau_m^\ell\times\tilde\tau_{m'}^\ell$ is contained on a set $H_{m,m',\ell}=\{(a,b)\in\R^N\times\R;\;
a=(m+m')2^{-\ell}+v,\;v=(v_1,v_2,\cdots,v_N),\; 0\le v_i\le2^{-\ell+1},\; 2^{-2\ell+2M}\le -|a|^2-\frac b\pi\le
3N2^{-2\ell+2M}\}.$ Consider the paraboloid defined by $-|a|^2-\frac b\pi=2^{-2\ell+2M}.$ Take $\Pi_{m,m',\ell}$ to be the tangent hyperplane to this paraboloid at the point of coordinates $(a_0,b_0),$ with $a_0= (m+m') 2^{-\ell},$ $b_0=-\pi |a_0|^2-2^{-2\ell+2M}$ (and passing through that point). Consider also the point $(a_1,b_1)$ with
$a_1=a_0+(2^{-\ell+1},2^{-\ell+1},\ldots,2^{-\ell+1})$ and $b_1=-\pi|a_1|^2-3 N 2^{-2\ell+2M}.$ Then, the rectangle
$R_{m,m',\ell}$ is defined as the only rectangle having a face contained in that hyperplane and the points $(a_0,b_0),$ and $(a_1,b_1)$ as opposite vertices. Due to the convexity of paraboloids, it follows that $H_{m,m',\ell}\subset R_{m,m',\ell}$ (see Figure~\ref{rectangle}). Moreover, one can also see that, for
small $c=c(N),$ $(1+c)R_{m,m',\ell}\subset \{ (a,b);\; a=(m+m')2^{-\ell}+v,\;v=(v_1,v_2,\ldots,v_N),\; |v_i|\le C(N)
2^{-\ell+1},\; C'(N) 2^{-2\ell+2M}\le-|a|^2-\frac b\pi\le C''(N)2^{-2\ell+2M}\}.$ Therefore, we have
$\sum_\vD\chi_{(1+c)R_{m,m',\ell}}\le C(N).$ Hence \eqref{stric1-1} in the case $\wh g\in X_{p,q}.$ Now, assume $g\in X_{p,q}.$ By density, it is sufficient to prove \eqref{stric1-1} for $g\in L^2(\R^N).$ By a straightforward calculation and the above result, we obtain that
$\|\vT(\:.\:)g\|_{L^q(\R^{N+1})}=\|\vT(\:.\:)\left(\vF^{-1}\ovl g\right)\|_{L^q(\R^{N+1})}\le C(N,p)\|g\|_{X_{p,q}}.$ Hence \eqref{stric1-1}.
\medskip
\end{vproof}

\begin{vproof}{of Theorem~\ref{stric2}.}
Notice first, that the second inequality follows from H{\"o}lder's. By homogeneity, we can assume that
$\|f\|_{L^2(\R^N)}=1.$ Then, it suffices to show that for any function $f\in L^2(\R^N)$ such that $\|f\|_{L^2(\R^N)}=1,$
$$
\sum_j \sum_k 2^{j\frac N2\frac{2-p}p q}\left(\int_{\tau_k^j}|f|^p\right)^{\frac{q}{p}}\le
C(p,q)\left[\sup_{j,k}\left\{2^{j\frac N2\frac{2-p}p}\left(\int_{\tau_k^j}|f|^p\right)^{\frac1p}\right\}\right]^{\alpha},
$$
where $\alpha=\mu pq$ and where $\mu$ has to be determined. Take $\alpha$ and $\beta$ such that $\frac2q<\beta<1,$ $\beta>\frac p2$ and $\alpha+q\beta=q.$ Then,
\begin{gather*}
\sum_j\sum_k 2^{j\frac N2\frac{2-p}p q}\left(\int_{\tau_k^j}|f|^p\right)^{\frac{q}{p}}
\le\left\{\sum_j\sum_k 2^{j\frac N2\frac{2-p}p\beta q}\left(\int_{\tau_k^j}|f|^p\right)^{\beta\frac{q}{p}}\right\}
\sup_{j,k}\left[2^{j\frac N2\frac{2-p}p}\left(\int_{\tau_k^j}|f|^p\right)^\frac{1}{p}\right]^\alpha.
\end{gather*}
We set $\mu=\frac{\alpha}{pq}=\frac{1-\beta}{p}\in\left(0,\frac{1}{p}\right).$ Hence, it is enough to show
\begin{gather*}
\sum_j\sum_k 2^{j\frac N2\frac{2-p}p\beta q}\left(\int_{\tau_k^j}|f|^p\right)^{\beta\frac{q}{p}}\le C(p,q).
\end{gather*}
We split the sum,
\begin{align*}
     & \; \sum_j \sum_k2^{j\frac N2\frac{2-p}p\beta q}\left(\int_{\tau_k^j}|f|^p\right)^{\beta\frac{q}{p}}
            \medskip \\
\le & \; C\sum_j\sum_k 2^{j\frac N2\frac{2-p}p\beta q}
            \left(\int_{\tau_k^j\cap\{|f|>2^{jN/2}\}}|f|^p\right)^{\beta\frac{q}{p}} \medskip \\
 +  & \; C\sum_j\sum_k 2^{j\frac N2\frac{2-p}p\beta q}
            \left(\int_{\tau_k^j\cap\{|f|\le 2^{jN/2}\}}|f|^p\right)^{\beta\frac{q}{p}}\stackrel{{\rm not}}{=}C(A+B),
\end{align*}
where $C=C(p,q).$ We study the first term. Set for each $j\in\Z,$ $f^j=f\chi_{\{|f|>2^{jN/2}\}}.$ Then,
\begin{gather*}
A=\sum_j\sum_k\left(2^{j\frac N2(2-p)}\int_{\tau_k^j}|f^j|^p\right)^{\beta\frac{q}{p}}.
\end{gather*}
Since $\beta q>2,$ we also have $\beta\frac{q}{p}>1.$ Then,
\begin{align*}
A\le & \;\left(\sum_j\sum_k2^{j\frac
N2(2-p)}\int_{\tau_k^j}|f^j|^p\right)^{\beta\frac{q}{p}}
            =\left(\sum_j2^{j\frac{N}{2}(2-p)}\int_{\R^N}|f^j|^p\right)^{\beta\frac{q}{p}} \medskip \\
   \le & \;\left(\int_{\R^N}|f|^p\sum_{\{j;\;|f|>2^{jN/2}\}}2^{j\frac N2(2-p)}\right)^{\beta\frac{q}{p}}.
\end{align*}
Since $2-p>0,$ we can sum the series and obtain
$$
A\le
C\left(\int_{\R^N}|f|^p|f|^{(2-p)}\right)^{\beta\frac{q}{p}}\le C
\left(\int_{\R^N}|f|^2\right)^{\beta\frac{q}{p}}\le C,
$$
by our assumption that $\|f\|_{L^2}=1.$ We now estimate $B.$ Set for any $j\in\Z,$ $f_j=f\chi_{\{|f|\le 2^{jN/2}\}}.$ Then,
$$
B=\sum_j\sum_k2^{j\frac N2\frac{2-p}p\beta q}\left(\int_{\tau_k^j}|f_j|^p\right)^{\beta\frac{q}{p}}
$$
We use H{\"o}lder's inequality with exponents $\frac{\beta q}{p}$ and $\frac{\beta q}{\beta q-p}.$ We obtain,
\begin{align*}
B\le & \; \sum_j\sum_k2^{j\frac N2\frac{2-p}p\beta
q}\int_{\tau_k^j}|f_j|^{\beta q}
               \left(|\tau_k^j|^\frac{\beta q-p}{\beta q}\right)^{\beta\frac{q}{p}} \medskip \\
    =  & \;\sum_j\sum_k2^{j\frac N2\frac{2-p}p\beta q}\int_{\tau_k^j}|f_j|^{\beta q}
               \left(2^{-jN\frac{\beta q-p}{\beta q}}\right)^{\beta\frac{q}{p}} \medskip \\
    =  & \; \sum_j\sum_k2^{jN(1-\beta \frac q2)}\int_{\tau_k^j}|f_j|^{\beta q}=\sum_j2^{jN(1-\beta \frac q2)}
               \int_{\R^N}|f_j|^{\beta q} \medskip \\
   \le & \int_{\R^N}|f|^{\beta q}\sum_{\{j;\;|f|\le 2^{jN/2}\}}2^{jN(1-\beta \frac q2)}.
\end{align*}
Since $1-\beta \frac q2<0,$ we sum the series to obtain
$$
B\le C\int_{\R^N}|f|^{\beta q}|f|^{(2-\beta q)}\le C
\int_{\R^N}|f|^2\le C,
$$
since $\|f\|_{L^2}=1.$
\\
We give an example to show that $L^2(\R^N)\neq X_{p,q}.$ Let
$$
f(x)=\frac1{|x|^{\frac N2} |\ln|x||^{\frac{1}{2}}}\chi_{\left(0,\frac12\right)^N}.$$
Then for any $1\le p<2$ and any $q>2,$ $f\in X_{p,q}$ but $f\not\in L^2(\R^N).$
\medskip
\end{vproof}

\section{Preliminary results}
\label{preli}

In this and next section, we follow Bourgain's arguments (\cite{MR99f:35184}). We have to modify them in the proof of Lemma~\ref{lemg}, because the Strichartz's exponent is not, in general, a natural number.

\begin{lem}
\label{lemf}
Let $f\in L^2(\R^N)\setminus\{0\}.$ Then for any $\eps>0,$ such that $\|\mathcal T(\cdot)f\|_{L^{\frac{2(N+2)}N}(\R\times\R^N)}\ge\eps,$ there exist $N_0\in\N$ with $N_0\le C(\|f\|_{L^2},N,\eps),$ $(A_n)_{1\le n\le N_0}\subset(0,\infty)$ and $(f_n)_{1\le n\le N_0}\subset L^2(\R^N)$ satisfying the following properties.
\begin{enumerate}
\item
\label{f1}
$\forall n\in[\![1,N_0]\!],$ $\supp\wh{f_n}\subset\tau_n,$ where $\tau_n\in\vC$ with $\ell(\tau_n)\le
C\|f\|_{L^2(\R^N)}^c\eps^{-\nu}A_n,$ and where the constants $C,$ $c$ and $\nu$ are positive and depend only on $N.$
\item
\label{f2}
$\forall n\in[\![1,N_0]\!],$ $|\wh{f_n}|<A_n^{-\frac{N}{2}}.$
\item
\label{f3}
$\|\vT(\: . \:)f-\dsp\vsum_{n=1}^{N_0}\vT(\: .\:)f_n\|_{L^\frac{2(N+2)}{N}(\R\times\R^N)}<\eps.$
\item
\label{f4}
$\|f\|_{L^2(\R^N)}^2=\sum\limits_{n=1}^{N_0}\|f_n\|_{L^2(\R^N)}^2+\|f-\sum\limits_{n=1}^{N_0}f_n\|_{L^2(\R^N)}^2.$
\end{enumerate}
\end{lem}

\noindent
The proof relies on the following lemma.

\begin{lem}
\label{lemlemf}
Let $g\in L^2(\R^N)$ and let $\eps>0$ be such that $\|\vT(\:.\:)g\|_{L^\frac{2(N+2)}{N}(\R\times\R^N)}\ge\eps.$ Then
there exist $h\in L^2(\R^N)$ and $A>0$ satisfying the following properties.
\begin{enumerate}
\item
\label{ff1}
$\supp\wh h\subset\tau,$ where $\tau\in\vC$ with $\ell(\tau)\le C\|g\|_{L^2(\R^N)}^c\eps^{-\nu}A,$ and where the constants $C,$ $c$ and $\nu$ depend only on $N.$
\item
\label{ff2}
$|\wh{h}|\le A^{-\frac{N}{2}}$ and $\|h\|_{L^2(\R^N)}^2\ge C\|g\|_{L^2(\R^N)}^{-a}\eps^b,$ where the constants $C,$ $a$ and $b$ depend only on $N.$
\item
\label{ff3}
$\|g-h\|_{L^2(\R^N)}^2=\|g\|_{L^2(\R^N)}^2-\|h\|_{L^2(\R^N)}^2.$
\end{enumerate}
\medskip
\end{lem}

\begin{proof*}
We distinguish 3 cases.
\medskip
\\
{\bf Case 1.} $\supp\wh g\subset [-1,1]^N.$ Then the function $h$
will also satisfy $\supp\wh h\subset\tau\subset [-1,1]^N.$
\\
Let $\eps>0$ and let $g$ be as in Lemma~\ref{lemlemf} such that $\supp\wh g\subset [-1,1]^N.$ It follows from Theorem~\ref{stric3} that
$$
\eps\le\|\vT(\:.\:)g\|_{L^\frac{2(N+2)}{N}(\R\times\R^N)}\le C\|g\|_{L^2(\R^N)}^{1-\mu p}\left[\vsup_{(j,k)\in\Z\times\Z^N}
2^{j\frac{N}{2}(2-p)}\int_{\tau_k^j}|\wh{g}(\xi)|^pd\xi\right]^\mu.
$$
So there exist $j\in\Z$ and $\tau\in\vC,$ with $\tau\subset[-1,1]^N$ and $\ell(\tau)=2^{-j},$ such that
\begin{gather}
\label{demlemlem1}
\vint_\tau|\wh{g}(\xi)|^pd\xi\ge C(\|g\|_{L^2(\R^N)}^{\mu p-1}\eps)^\frac{1}{\mu}2^{-j\frac{N}{2}(2-p)}.
\end{gather}
Let $M=\left((C\|g\|_{L^2(\R^N)}^{\mu(p-2)-1}\eps)^\frac{1}{\mu}2^{-j\frac{N}{2}(2-p)-1}\right)^\frac{1}{p-2},$
where $C$ is the constant in (\ref{demlemlem1}). Then by Plancherel's Theorem,
\begin{gather}
\label{demlemlem2}
\vint_{\tau\cap\{|\wh g|\ge M\}}|\wh g(\xi)|^pd\xi=M^{p-2}\vint_{\tau\cap\{|\wh g|\ge M\}}|\wh g(\xi)|^pM^{2-p}d\xi
\le M^{p-2}\int |\wh g|^p|\wh g|^{2-p} = M^{p-2}\|g\|_{L^2(\R^N)}^2.
\end{gather}
It follows from (\ref{demlemlem1})--(\ref{demlemlem2}) that
\begin{align*}
& \vint_{\tau\cap\{|\wh g|<M\}}|\wh g(\xi)|^pd\xi
    =\vint_\tau|\wh g(\xi)|^pd\xi-\vint_{\tau\cap\{|\wh g|\ge M\}}|\wh g(\xi)|^pd\xi \medskip \\
& \ge (C\|g\|_{L^2(\R^N)}^{\mu p-1}\eps)^\frac{1}{\mu}2^{-j\frac{N}{2}(2-p)}-M^{p-2}\|g\|_{L^2(\R^N)}^2 \medskip \\
& \ge C\eps^\frac{1}{\mu}2^{-j\frac{N}{2}(2-p)}\|g\|_{L^2(\R^N)}^{-\frac{1-\mu p}{\mu}}.
\end{align*}
By H{\"o}lder's inequality and the above estimate, we get
$$
C\eps^\frac{1}{\mu}2^{-j\frac{N}{2}(2-p)}\|g\|_{L^2(\R^N)}^{-\frac{1-\mu p}{\mu}}
\le\vint_{\tau\cap\{|\wh g|<M\}}|\wh g(\xi)|^pd\xi\le\left(\:\vint_{\tau\cap\{|\wh g|<M\}}
|\wh g(\xi)|^2d\xi\right)^\frac{p}{2}|\tau|^\frac{2-p}{2}.
$$
Since $|\tau|=2^{-jN},$ we then obtain,
\begin{gather}
\label{demlemlem3}
\vint_{\tau\cap\{|\wh g|<M\}}|\wh g(\xi)|^2d\xi\ge C\|g\|_{L^2(\R^N)}^{-\frac{2(1-\mu p)}{\mu p}}\eps^\frac{2}{\mu p}.
\end{gather}
Let $h\in L^2(\R^N)$ be such that $\wh h=\wh g\chi_{\tau\cap\{|\wh g|<M\}}$ and let $A=M^{-\frac{2}{N}}.$ Then $\supp\wh h\subset\tau\subset [-1,1]^N$ with $\ell(\tau)=2^{-j}=C\|g\|_{L^2(\R^N)}^\frac{2\mu(2-p)+2}{N\mu(2-p)}\eps^{-\frac{2}{N\mu(2-p)}}A.$ So we have~\ref{ff1}, and~\ref{ff2} follows from (\ref{demlemlem3}). Since $\wh h$ and $\wh g-\wh h$ have disjoint supports, \ref{ff3} follows.
\medskip
\\
{\bf Case 2.} $\supp\wh g\subset [-M,M]^N$ for some $M>0.$ Then $h$ will also satisfy $\supp\wh h\subset\tau\subset [-M,M]^N.$
\\
Let $\eps>0$ and let $g$ be as in the Lemma~\ref{lemlemf} such that $\supp\wh g\subset [-M,M]^N$ for some $M>0.$ Let $g'\in L^2(\R^N)$ be such that $\wh g'(\xi)=M^\frac{N}{2}\wh g(M\xi).$ Then
$\supp\wh g'\subset [-1,1]^N$ and so we may apply the Case~1 to $g'.$ Thus there exist $h'\in L^2(\R^N),$ $\tau'\in\vC$ and $A'>0$ satisfying \ref{ff1}--\ref{ff3}. We define $h\in L^2(\R^N)$ by $\wh h(\xi)=M^{-\frac{N}{2}}\wh h'\left(\frac{\xi}{M}\right).$ Then $\|g\|_{L^2(\R^N)}=\|g'\|_{L^2(\R^N)}$ and $\|h\|_{L^2(\R^N)}=\|h'\|_{L^2(\R^N)}.$ In particular, second part of~\ref{ff2} holds for $g$ and $h.$ Setting $\tau=M\tau',$ it follows that $\supp\wh h\subset\tau\subset [-M,M]^N$ and $\ell(\tau)=M\ell(\tau')\le C\|g\|_{L^2(\R^N)}^q\eps^\nu MA'.$ So $h$ satisfies~\ref{ff1} with $A=MA'.$ Finally, $|\wh h|<M^{-\frac{N}{2}}A'^{-\frac{N}{2}}=A^{-\frac{N}{2}},$ which
implies~\ref{ff2}. Finally,~\ref{ff3} follows from the similar identity for $\wh g'$ and $\wh h'.$
\medskip
\\
{\bf Case 3.} General case. \\
Let $\eps>0$ and let $g$ be as in the Lemma~\ref{lemlemf}. For $M>0,$ we define $u_M\in L^2(\R^N)$ by $\wh{u_M}=\wh g\chi_{[-M,M]^N}.$ It follows from Strichartz's estimate (\ref{stri1}) and Plancherel's Theorem that
$$
\|\vT(\:.\:)(u_M-g)\|_{L^\frac{2(N+2)}{N}(\R\times\R^N)}\le C\|u_M-g\|_{L^2(\R^N)}
=C\|\wh{u_M}-\wh g\|_{L^2(\R^N)}\xrightarrow{M\tends\infty}0.
$$
Then there exists $M_0>0$ such that
$$
\|\vT(\:.\:)u_{M_0}\|_{L^\frac{2(N+2)}{N}(\R\times\R^N)}\ge\dfrac{\eps}{2}.
$$
Setting $g_0=u_{M_0},$ we apply the Case~2 to $g_0,$ obtaining $h.$ Since $\|g_0\|_{L^2(\R^N)}\le\|g\|_{L^2(\R^N)},$ Properties~\ref{ff1} and \ref{ff2} are clear for $g$ and $h.$ Also, Property~\ref{ff3} holds for $g$ and $h,$ again because the disjointness of supports. This achieves the proof of the lemma.
\medskip
\end{proof*}

\begin{vproof}{of Lemma~\ref{lemf}.}
Let $f\in L^2(\R^N)\setminus\{0\}$ and let $\eps>0$ be such that
$$
\|\vT(\:.\:)f\|_{L^\frac{2(N+2)}{N}(\R\times\R^N)}\ge\eps.
$$
We apply Lemma~\ref{lemlemf} to $f.$ Let $h\in L^2(\R^N),$ $\tau\in\vC,$ $A>0,$ $a=a(N)>0,$ $b=b(N)>0,$ $c=c(N)>0$ and $\nu=\nu(N)>0$ be given by Lemma~\ref{lemlemf}. We set $f_1=h,$ $\tau_1=\tau$ and $A_1=A.$ By Lemma~\ref{lemlemf}, we have
\begin{gather}
\label{demlemf1}
\ell(\tau_1)\le C\|f\|_{L^2}^c\eps^{-\nu}A_1, \\
\label{demlemf2}
\|f-f_1\|_{L^2}^2=\|f\|_{L^2}^2-\|f_1\|_{L^2}^2,
\quad
\|f-f_1\|_{L^2}^{-a}\ge\|f\|_{L^2}^{-a}
\quad\mbox{and}\quad
\|f_1\|_{L^2}^2\ge C\|f\|_{L^2}^{-a}\eps^b.
\end{gather}
Now, we may assume that
$$
\|\vT(\:.\:)f-\vT(\:.\:)f_1\|_{L^\frac{2(N+2)}{N}(\R\times\R^N)}\ge\eps,
$$
otherwise we set $N_0=1$ and the proof is finished. So we may apply Lemma~\ref{lemlemf} to $g=f-f_1.$ Let $h\in L^2(\R^N),$ let $\tau\in\vC$ and let $A>0$ be given by Lemma~\ref{lemlemf}. We set $f_2=h,$ $\tau_2=\tau$ and $A_2=A.$ By Lemma~\ref{lemlemf} and (\ref{demlemf2}), we have
\begin{gather}
\label{demlemf3}
\ell(\tau_2)\le C\|f-f_1\|_{L^2}^c\eps^{-\nu}A_2\le C\|f\|_{L^2}^c\eps^{-\nu}A_2, \\
\label{demlemf4}
\|f-(f_1+f_2)\|_{L^2}^2=\|f-f_1\|_{L^2}^2-\|f_2\|_{L^2}^2=\|f\|_{L^2}^2-(\|f_1\|_{L^2}^2+\|f_2\|_{L^2}^2), \\
\label{demlemf5}
\|f_2\|_{L^2}^2\ge C\|f-f_1\|_{L^2}^{-a}\eps^b\ge C\|f\|_{L^2}^{-a}\eps^b.
\end{gather}
We repeat the process as long as
$$
\|\vT(\: . \:)f-\dsp\vsum_{j=1}^{k-1}\vT(\: . \:)f_j\|_{L^\frac{2(N+2)}{N}(\R\times\R^N)}\ge\eps,
$$
applying Lemma~\ref{lemlemf}  to $g=f-\dsp\vsum_{j=1}^{k-1}f_j.$ Then, by (\ref{demlemf1})--(\ref{demlemf5}), we obtain functions $f_1,\ldots,f_n$ satisfying Properties~\ref{f1} and \ref{f2} of
Lemma~\ref{lemf} and
\begin{gather}
\label{demlemf6}
\|f-\vsum_{j=1}^kf_j\|_{L^2}^2=\|f\|_{L^2}^2-\vsum_{j=1}^k\|f_j\|_{L^2}^2, \medskip \\
\label{demlemf7}
\|f_k\|_{L^2}^2\ge C\|f\|_{L^2}^{-a}\eps^b,
\end{gather}
for any $ k\in[\![1,n]\!],$ for some $n\ge2.$ From Strichartz's estimate (\ref{stri1}) and (\ref{demlemf6})--(\ref{demlemf7}), we obtain
\begin{align*}
     & \;\|\vT(\: . \:)f-\dsp\vsum_{j=1}^n\vT(\: . \:)f_j\|_{L^\frac{2(N+2)}{N}(\R\times\R^N)}^2 \medskip \\
\le & \;C\|f-\vsum_{j=1}^nf_j\|_{L^2}^2\le C(\|f\|_{L^2}^2-Cn\|f\|_{L^2}^{-a}\eps^b)\xrightarrow{n\tends\infty}
         -\infty.
\end{align*}
So the process stops for some $n\le C(\|f\|_{L^2},N,\eps).$ We set $N_0=n$ and the proof is achieved.
\medskip
\end{vproof}

\begin{lem}
\label{lemg}
Let $g\in L^2(\R^N),$ let $\tau\in\vC,$ let $A>0$ and let $C_0>0$ be such that $\supp\wh g\subset\tau,$ $\ell(\tau)\le C_0A$ and $|\wh g|<A^{-\frac{N}{2}}.$ Let $\xi_0$ be the center of $\tau.$ Then for any $\eps>0,$ there exist $N_1\in\N$ with $N_1\le C(N,C_0,\eps)$ and $(Q_n)_{1\le n\le N_1}\subset\R\times\R^N$ with
\begin{gather}
\label{lemg1}
Q_n=\left\{(t,x)\in\R\times\R^N;\; t\in I_n \mbox{ and } (x-4\pi t\xi_0)\in C_n \right\},
\end{gather}
where $I_n\subset\R$ is an interval with $|I_n|=\dfrac{1}{A^2}$ and $C_n\in\vC$ with $\ell(C_n)=\dfrac{1}{A}$ such that
$$
\bigg(\:\vint_{\R^{N+1}\setminus\vcup_{n=1}^{N_1}Q_n}|(\vT(t))g(x)|^\frac{2(N+2)}{N}dtdx\bigg)^\frac{N}{2(N+2)}
<\eps.
$$
\end{lem}

Notice that the functions $f_n$ obtained in Lemma~\ref{lemf} satisfy the hypothesis of Lemma~\ref{lemg}.
\bigskip

\begin{vproof}{of Lemma \ref{lemg}.}
We define $g'\in L^2(\R^N)$ by $\wh{g'}(\xi')=A^\frac{N}{2}\wh g(\xi_0+A\xi').$ Then $\|g'\|_{L^2}=\|g\|_{L^2},$ $|\wh{g'}|<1$ and $\supp\wh{g'}\subset\left[-\frac{C_0}{2},\frac{C_0}{2}\right]^N.$ It follows from
(\ref{demlemg1}) applied to $g'$ that
\begin{align*}
   & |(\vT(A^2t)g')(A(x-4\pi t\xi_0))|
       =\bigg|\:\vint_{\left(-\frac{C_0}{2},\frac{C_0}{2}\right)^N}e^{2i\pi\left(A(x-4\pi t\xi_0).\xi-2\pi A^2t|\xi|^2\right)}
       \wh{g'}(\xi)d\xi\bigg| \medskip \\
= &
\;A^\frac{N}{2}\bigg|\:\vint_{\left(-\frac{C_0}{2},\frac{C_0}{2}\right)^N}e^{2i\pi
       \left(A(x-4\pi t\xi_0).\xi-2\pi A^2t|\xi|^2\right)}\wh g(\xi_0+A\xi)d\xi\bigg| \medskip \\
= & \;A^{-\frac{N}{2}}|(\vT(t)g)(x)|,
\end{align*}
where the last identity follows from the change of variables
$\zeta=\xi_0+A\xi.$ Setting
\begin{gather}
\begin{cases}
\label{deftx}
t'=A^2t, \\
x'=A(x-4\pi t\xi_0),
\end{cases}
\end{gather}
we then have
\begin{gather}
\label{demlemg3}
|(\vT(t)g)(x)|=A^\frac{N}{2}|(\vT(t')g')(x')|.
\end{gather}
By (\ref{demlemg1}),
\begin{gather}
\label{demlemg4}
|(\vT(t)g')(x)|
=\bigg|\:\vint_{\left(-\frac{C_0}{2},\frac{C_0}{2}\right)^N}\wh{g'}(\zeta)e^{2i\pi(x.\zeta-2\pi
t|\zeta|^2)}d\zeta\bigg|.
\end{gather}
By~(\ref{adjrest}) (with $g'$ in the place of $g)$ and Corollary~1.2 of~Tao~\cite{MR2033842}, we obtain
\begin{equation}
   \label{demlemg5}
        \|\vT(\:.\:)g'\|_{L^q(\R\times\R^N)}\le C(N,q)\|\wh{g'}\|_{L^p(\R^N)}
            =C(N,q)\|\wh{g'}\|_{L^p\left(\left(-\frac{C_0}{2},\frac{C_0}{2}\right)^N\right)},
\end{equation}
for any $q>\frac{2(N+3)}{(N+1)}$ and any $p\ge1$ such that $q=\frac{N+2}{N}p'.$ Let $p'=p'(N)\in(1,2)$ be such that
$$
\dfrac{2(N+3)}{(N+1)}<\dfrac{N+2}{N}p'<\dfrac{2(N+2)}{N}.
$$
Thus $q=q(N)=\dfrac{N+2}{N}p'<\dfrac{2(N+2)}{N}$ and it follows from (\ref{demlemg5}) that and H{\"o}lder's inequality that
\begin{gather*}
\|\vT(\:.\:)g'\|_{L^q(\R\times\R^N)}\le C(N)\|\wh{g'}\|_{L^p\left(\left(-\frac{C_0}{2},\frac{C_0}{2}\right)^N\right)}
\le C(N)\bigg|\bigg(-\frac{C_0}{2},\frac{C_0}{2}\bigg)^N\bigg|^\frac{1}{p}
\|\wh{g'}\|_{L^\infty\left(\left(-\frac{C_0}{2},\frac{C_0}{2}\right)^N\right)},
\end{gather*}
so that
$$
\|\vT(\:.\:)g'\|_{L^q(\R\times\R^N)}\le C(C_0,N).
$$
This estimate implies that for any $\lambda>0,$
\begin{align*}
   & \;\vint_{\{|\vT(\:.\:)g'|<\lambda\}}|\vT(t')g'(x')|^\frac{2(N+2)}{N}dt'dx' \medskip \\
= &
\;\vint_{\{|\vT(\:.\:)g'|<\lambda\}}|\vT(t')g'(x')|^{\left(\frac{2(N+2)}{N}-q\right)+q}dt'dx'
        \le C(C_0,N)\lambda^{\frac{2(N+2)}{N}-q}.
\end{align*}
So there exists $\lambda_0=\lambda_0(N,C_0,\eps)\in(0,1)$ small enough such that
\begin{gather}
\label{demlemg7}
\vint_{\{|\vT(\:.\:)g'|<2\lambda_0\}}|\vT(t')g'(x')|^{\frac{2(N+2)}{N}}dt'dx'<\eps^\frac{2(N+2)}{N},
\end{gather}
Since
$\supp\wh{g'}\subset\left[-\frac{C_0}{2},\frac{C_0}{2}\right]^N$ and $\|\wh g'\|_{L^\infty}\le1,$ it follows from formula
(\ref{demlemg1}) that for any $(t',x')\in\R\times\R^N$ and any $(t'',x'')\in\R\times\R^N,$
$$
|\vT(t')g'(x')-\vT(t'')g'(x'')|\le C(|t'-t''|+|x'-x''|),
$$
where $C=C(C_0,N)\ge1.$ So for such a constant, if
$(t',x')\in\{|\vT(\:.\:)g'|\ge2\lambda_0\}$ and if $(t'',x'')\in\R\times\R^N$ is such that
$|t'-t''|\le\frac{\lambda_0}{2C}<\frac{1}{2}$ and $|x'-x''|\le\frac{\lambda_0}{2C} <\frac{1}{2}$ then
$|\vT(t'')g(x'')|\ge\lambda_0,$ that is $(t'',x'')\in\{|\vT(\:.\:)g'|\ge\lambda_0\}.$ So there exist a set $R$ and a family $(P_r)_{r\in R}=(J_r,K_r)_{r\in R}\subset\R\times\R^N,$ where $J_r\subset\R$ is a closed interval of center $t'\in\R$ with $|J_r|=\frac{\lambda_0}{C}$ and $K_r\in\vC$ of center $x'\in\R^N$ with $\ell(K_r)=\frac{\lambda_0}{C}$ and $(t',x')\in\{|\vT(\:.\:)g'|\ge2\lambda_0\},$ such that
\begin{gather}
\label{demlemg8}
\forall(r,s)\in R\times R \mbox{ with } r\neq s, \; {\rm Int }(P_r)\cap{\rm Int}(P_s)=\emptyset, \\
\label{demlemg9}
\{|\vT(\:.\:)g'|\ge2\lambda_0\}\subset\bigcup_{r\in R}P_r\subset\{|\vT(\:.\:)g'|\ge\lambda_0\},
\end{gather}
where ${\rm Int }(P_r)$ denotes the interior of the set $P_r.$ We set $N_1=\# R.$ It follows from
(\ref{demlemg8})--(\ref{demlemg9}) and Strichartz's estimate (\ref{stri1}) that,
\begin{align*}
     & \;N_1\left(\frac{\lambda_0}{C}\right)^{N+1}=\left|\bigcup_{r\in R}P_r\right|
         \le\left|\left\{|(\vT(\:.\:)g')|\ge\lambda_0\right\}\right| \medskip \\
\le & \;\lambda_0^{-\frac{2(N+2)}{N}}\|\vT(\:.\:)g'\|_{L^\frac{2(N+2)}{N}(\R\times\R^N)}^\frac{2(N+2)}{N}
         \le C\lambda_0^{-\frac{2(N+2)}{N}}\|g\|_{L^2}^\frac{2(N+2)}{N},
\end{align*}
from which we deduce that $N_1<\infty$ and $N_1\le C(\|g\|_{L^2},N,C_0,\eps).$ Actually, since our hypothesis implies that $\|g\|_{L^2}\le C_0^{N/2},$ we can write also $N_1\le C(N,C_0,\eps).$ For any $n\in[\![1,N_1]\!],$ let $(t_n,x_n)$ be the center of $P_n,$ let $I_n\subset\R$ be the interval of center $\frac{t_n}{A^2}$ with $|I_n|=\frac{1}{A^2},$ let $I'_n=A^2I_n,$ let $C_n\in\vC$ of center $\frac{1}{A}x_n$ with $\ell(C_n)=\frac{1}{A},$ let $C'_n=AC_n$ and let $Q_n$ be defined by~(\ref{lemg1}). Then $\bigcup\limits_{n=1}^{N_1}P_n\subset\bigcup\limits_{n=1}^{N_1}(I'_n\times C'_n),$ which yields with (\ref{demlemg7}) and (\ref{demlemg9}),
\begin{gather}
\label{demlemg10}
\vint_{\R^{N+1}\setminus\bigcup\limits_{n=1}^{N_1}(I'_n\times C'_n)}
|\vT(t')g'(x')|^\frac{2(N+2)}{N}dt'dx'<\eps^\frac{2(N+2)}{N}.
\end{gather}
By (\ref{demlemg3}),
\begin{gather*}
\vint_{\R^{N+1}\setminus\bigcup\limits_{n=1}^{N_1}Q_n}|\vT(t)g(x)|^\frac{2(N+2)}{N}dtdx
=A^{N+2}\vint_{\R^{N+1}\setminus\bigcup\limits_{n=1}^{N_1}Q_n}|\vT(t')g'(x')|^\frac{2(N+2)}{N}dtdx
\end{gather*}
But $(t,x)\in Q_n\iff(t',x')\in I'_n\times C'_n,$ and so we deduce from the above estimate and (\ref{deftx}) that
\begin{gather}
\label{demlemg11}
\vint_{\R^{N+1}\setminus\bigcup\limits_{n=1}^{N_1}Q_n}|\vT(t)g(x)|^\frac{2(N+2)}{N}dtdx
=\vint_{\R^{N+1}\setminus\bigcup\limits_{n=1}^{N_1}(I'_n\times
C'_n)}|\vT(t')g(x')|^\frac{2(N+2)}{N}dt'dx'.
\end{gather}
Putting together (\ref{demlemg10}) and (\ref{demlemg11}), we obtain the desired result.
\medskip
\end{vproof}

\section{Mass concentration}
\label{mass}

\begin{prop}
\label{propmass}
Let $\gamma\in\R\setminus\{0\},$ let $u_0\in L^2(\R^N)\setminus\{0\}$ and let
$$
u\in C((-T_\m,T_\M);L^2(\R^N))\cap L_\loc^\frac{2(N+2)}{N}((-T_\m,T_\M);L^\frac{2(N+2)}{N}(\R^N))
$$
be the maximal solution of $(\ref{nls})$ such that $u(0)=u_0.$ Then there exists $\eta_0=\eta_0(N,|\gamma|)>0$ satisfying the following properties. Let $(T_0,T_1)\subset(-T_\m,T_\M)$ be an interval and let
\begin{gather}
\label{propmass1}
\eta=\|u\|_{L^\frac{2(N+2)}{N}((T_0,T_1)\times\R^N)}.
\end{gather}
If $\eta\in(0,\eta_0]$ then there exist $t_0\in(T_0,T_1)$ and $c\in\R^N$  such that
\begin{gather}
\label{propmass2}
\|u(t_0)\|_{L^2(B(c,R))}\ge\eps,
\end{gather}
where $R=\min\big\{(T_1-t_0)^\frac{1}{2}, (t_0-T_0)^{\frac12}\big\}$ and $\eps=\eps(\|u_0\|_{L^2},N,\eta)>0.$
\end{prop}

\begin{proof*}
Let $\gamma,$ $u_0,$ $u$ and $(T_0,T_1)$ be as in the Proposition~\ref{propmass}. Let $\eta>0$ be as in
(\ref{propmass1}). By (\ref{inteq}), we have
\begin{gather}
\label{dempropmass1}
\forall t\in(-T_\m,T_\M),\; u(t)=\vT(t-T_0)u(T_0)+i\gamma\int_{T_0}^t(\vT(t-s)\{|u|^\frac{4}{N}u\})(s)ds.
\end{gather}
Setting for any $t\in(-T_\m,T_\M),$ $\Phi_u(t)=i\gamma\int_{T_0}^t(\vT(t-s)\{|u|^\frac{4}{N}u\})(s)ds$ and applying Strichartz's estimate (\ref{stri2}), we get with (\ref{propmass1})
\begin{gather}
\label{dempropmass2}
\|\Phi_u\|_{L^\frac{2(N+2)}{N}((T_0,T_1)\times\R^N)}
\le C_1\|u\|_{L^\frac{2(N+2)}{N}((T_0,T_1)\times\R^N)}^\frac{N+4}{N}=C_1\eta^\frac{N+4}{N},
\end{gather}
where $C_1=C_1(N,|\gamma|)\ge1.$ For every $a,b\ge0,$ $(a+b)^\alpha\le C(\alpha)(a^\alpha+b^\alpha),$ where
$C(\alpha)=1$ if $0<\alpha\le1$ and $C(\alpha)=2^{\alpha-1}$ if $\alpha\ge1.$ Let $C_2$ be such a constant for
$\alpha=\frac{4}{N}.$ We choose $\eta_0=\eta_0(N,|\gamma|)>0$ small enough to have
\begin{gather}
\label{dempropmass4}
2(2C_1)^\frac{4}{N}C_2\eta_0^\frac{16}{N^2}\le 1.
\end{gather}
Assume that $\eta\le\eta_0.$   We proceed in 3 steps.
\medskip
\\
{\bf Step 1.} We show that, there exist $f_0\in L^2(\R^N),$ $A>0$ and $\tau\in\vC$ of center $\xi_0\in\R^N$ satisfying $\supp\wh{f_0}\subset\tau,$ $\ell(\tau)\le C(\|u_0\|_{L^2},N,\eta)A$ and $|\wh{f_0}|<A^{-\frac{N}{2}},$ and
there exist an interval $I\subset\R$ and $K\in\vC,$ with $|I|=\dfrac{1}{A^2}$ and $\ell(K)=\dfrac{1}{A},$ such that for
$Q\subset\R\times\R^N$ defined by
\begin{gather*}
Q=\left\{(t,x)\in\R\times\R^N;\; t\in I \mbox{ and } (x-4\pi t\xi_0)\in K\right\},
\end{gather*}
we have
\begin{gather}
\label{dempropmass3}
\iint\limits_{((T_0,T_1)\times\R^N)\cap Q}|u(t,x)|^2|\vT(t-T_0)f_0(x)|^\frac{4}{N}dtdx\ge C\eta^\frac{2(N+2)}{N},
\end{gather}
where $C=C(\|u_0\|_{L^2},N,\eta).$
\\
To prove this claim, we apply Lemma~\ref{lemf} to $f=u(T_0)$ with $\eps_0=\eta^\frac{N+4}{N}.$ Note that, by \eqref{propmass1}, \eqref{dempropmass1}, \eqref{dempropmass2}, \eqref{dempropmass4} and time translation, we have that
$$
\|\mathcal T(\cdot)u(T_0)\|_{L^{\frac{2(N+2)}N}(\R\times\R^N)}=\|\mathcal
T(\cdot-T_0)u(T_0)\|_{L^{\frac{2(N+2)}N}(\R\times\R^N)}\ge\eta/2\ge\eps_0.
$$
It follows from H{\"o}lder's inequality (with $p=\frac{N+2}{N}$ and $p'=\frac{N+2}{2}),$
(\ref{dempropmass1})--(\ref{dempropmass2}) and Lemma~\ref{lemf} that
\begin{align*}
     & \;\iint\limits_{T_0\:\R^N}^{\;\;T_1}|u(t,x)|^2\left|u(t,x)-\sum_{n=1}^{N_0}\vT(t-T_0)f_n(x)\right|^\frac{4}{N}dtdx
         \medskip \\
\le & \;\|u\|_{L^\frac{2(N+2)}{N}((T_0,T_1)\times\R^N)}^2
         \|u-\dsp\vsum_{n=1}^{N_0}\vT(\: \cdot-T_0 \:)f_n\|_{L^\frac{2(N+2)}{N}((T_0,T_1)\times\R^N)}^\frac{4}{N} \medskip \\
\le & \;\eta^2\left(\|\vT(\: . \:)u(T_0)-\dsp\vsum_{n=1}^{N_0}\vT(\: . \:)f_n\|_{L^\frac{2(N+2)}{N}(\R\times\R^N)}
         +C_1\|u\|_{L^\frac{2(N+2)}{N}((T_0,T_1)\times\R^N)}^\frac{N+4}{N}\right)^\frac{4}{N} \medskip \\
\le & \;C_1^\frac{4}{N}\eta^2(\eps_0+\eta^\frac{N+4}{N})^\frac{4}{N}
         \le(2C_1)^\frac{4}{N}\eta_0^\frac{16}{N^2}\eta^\frac{2(N+2)}{N}\le\frac{1}{2C_2}\eta^\frac{2(N+2)}{N}.
\end{align*}
The above estimate and (\ref{propmass1}) yield
\begin{align*}
\eta^\frac{2(N+2)}{N}
 = & \;\iint\limits_{T_0\:\R^N}^{\;\;T_1}|u(t,x)|^2\left|\left(u(t,x)-\sum_{n=1}^{N_0}\vT(t-T_0)f_n(x)\right)
        +\sum_{n=1}^{N_0}\vT(t-T_0)f_n(x)\right|^\frac{4}{N}dtdx \medskip \\
\le & \;C_2\left(\frac{1}{2C_2}\eta^\frac{2(N+2)}{N}
        +\iint\limits_{T_0\:\R^N}^{\;\;T_1}|u(t,x)|^2\left|\sum_{n=1}^{N_0}\vT(t-T_0)f_n(x)\right|^\frac{4}{N}dtdx\right),
\end{align*}
which gives
\begin{gather}
\label{dempropmass6}
\iint\limits_{T_0\:\R^N}^{\;\;T_1}|u(t,x)|^2\left|\sum_{n=1}^{N_0}\vT(t-T_0)f_n(x)\right|^\frac{4}{N}dtdx
\ge\frac{1}{2C_2}\eta^\frac{2(N+2)}{N}.
\end{gather}
By Lemma~\ref{lemf} and conservation of charge, $N_0\le C(\|u_0\|_{L^2},N,\eta).$ It follows from (\ref{dempropmass6}) that there exists $n_0\in[\![1,N_0]\!]$ such that
\begin{gather}
\label{dempropmass7}
\iint\limits_{T_0\:\R^N}^{\;\;T_1}|u(t,x)|^2\left|\vT(t-T_0)f_{n_0}(x)\right|^\frac{4}{N}dtdx\ge C\eta^\frac{2(N+2)}{N},
\end{gather}
where $C=C(\|u_0\|_{L^2},N,\eta).$ Set $A=A_{n_0},$ $\tau=\tau_{n_0}$ and $C_0=C(N)\|u_0\|_{L^2}^{c(N)}\eps_0^{-\nu(N)},$ where we have used the notations of Lemma~\ref{lemf}. Let $\xi_0\in\R^N$ be the center of $\tau_{n_0}.$ We apply Lemma~\ref{lemg} to $g=f_{n_0}$ and $\eps_1=\left(\frac{C}{2}\right)^\frac{N}{4}\eta,$ where $C$ is the constant in (\ref{dempropmass7}). It follows from
H{\"o}lder's inequality (with $p=\frac{N+2}{N}$ and $p'=\frac{N+2}{2}),$ (\ref{propmass1}) and Lemma~\ref{lemg} that
\begin{align*}
     & \;\iint\limits_{((T_0,T_1)\times\R^N)\setminus\vcup_{n=1}^{N_1}Q_n}|u(t,x)|^2
        \left|\vT(t-T_0)f_{n_0}(x)\right|^\frac{4}{N}dtdx \medskip \\
\le & \;\|u\|_{L^\frac{2(N+2)}{N}((T_0,T_1)\times\R^N)}^2
         \|\vT(\: . \:)f_{n_0}\|_{L^\frac{2(N+2)}{N}(\R^{N+1}\setminus\vcup_{n=1}^{N_1}Q_n)}^\frac{4}{N} \medskip \\
\le & \;\eta^2\eps_1^\frac{4}{N}=\frac{C}{2}\eta^\frac{2(N+2)}{N}.
\end{align*}
The above estimate with (\ref{dempropmass7}) yield
\begin{gather}
\label{dempropmass8}
\iint\limits_{((T_0,T_1)\times\R^N)\cap(\vcup_{n=1}^{N_1}Q_n)}|u(t,x)|^2
\left|\vT(t-T_0)f_{n_0}(x)\right|^\frac{4}{N}dtdx\ge C\eta^\frac{2(N+2)}{N},
\end{gather}
where $C=C(\|u_0\|_{L^2},N,\eta).$ By Lemma~\ref{lemg}, $N_1\le C(\|u_0\|_{L^2},N,\eta).$  With (\ref{dempropmass8}), this implies that there exists $n_1\in[\![1,N_1]\!]$ such that
\begin{gather}
\label{dempropmass13}
\iint\limits_{((T_0,T_1)\times\R^N)\cap Q_{n_1}}|u(t,x)|^2\left|\vT(t-T_0)f_{n_0}(x)\right|^\frac{4}{N}dtdx
\ge C\eta^\frac{2(N+2)}{N},
\end{gather}
where $C=C(\|u_0\|_{L^2},N,\eta).$ Hence we obtain the Step~1
claim with $f_0=f_{n_0},$ $I=I_{n_1},$ $K=C_{n_1}$ and
$Q=Q_{n_1}.$
\medskip
\\
{\bf Step 2.}
We show that $\dfrac{1}{A}\le C(T_1-T_0)^\frac{1}{2}$ and $\vsup_{t\in\R}\|\vT(t-T_0)f_0\|_{L^\infty(\R^N)}\le
CA^\frac{N}{2},$ where $C=C(\|u_0\|_{L^2},N,\eta).$
\\
By (\ref{demlemg1}) and Step~1, $|\vT(t-T_0)f_0|\le\dsp\int_\tau|\wh f_0(\xi)|d\xi\le A^{-\frac{N}{2}}\int_\tau1d\xi\le CA^\frac{N}{2},$ which yields second part of Step~2. Using this estimate, Step~1 and conservation of charge, we deduce
\begin{align*}
C\eta^\frac{2(N+2)}{N}\le & \;\iint\limits_{((T_0,T_1)\times\R^N)\cap Q}|u(t,x)|^2|\vT(t-T_0)f_0(x)|^\frac{4}{N}dxdt
                                                  \medskip \\
                                        \le & \;CA^2\iint\limits_{((T_0,T_1)\times\R^N)\cap Q}|u(t,x)|^2dxdxt
                                                 \le CA^2\iint\limits_{T_0\:\R^N}^{\;\;T_1}|u(t,x)|^2dxdt \medskip \\
                                       \le & \;CA^2\|u_0\|_{L^2}^2(T_1-T_0).
\end{align*}
Hence we obtain the Step~2 claim.
\medskip
\\
{\bf Step 3.} Conclusion. \\
Let $K\in \mathcal C,$ $I$ and $Q$ be as in Step~1, and let $\eta\prime=C\eta^\frac{2(N+2)}{N},$ where $C$ is the constant of \eqref{dempropmass13}. Let $K(t)=K+4\pi t\xi_0$ and let $\kappa>0$ be small enough to be chosen later. It follows from Step~1, Step~2 and H{\"o}lder's inequality (with $p=\frac{N+2}{N}$ and $p'=\frac{N+2}{2}),$ that
\begin{align*}
\eta\prime\le & \;\iint\limits_{((T_0,T_1)\times\R^N)\cap Q}|u(t,x)|^2\left|\vT(t-T_0)f_0(x)\right|^\frac{4}{N}dxdt
                           \medskip \\
                  \le & \; \|\vT(\cdot-T_0)f_0\|_{L^\infty}^\frac{4}{N}\vint_{I\cap\left(T_0,T_1\right)}
                            \left(\int_{K(t)}|u(t,x)|^2dx\right)dt
                            \medskip\\
                   \le & \; CA^2\vint_{I\cap\left(T_0,T_1\right)}
                             \left(\int_{K(t)}|u(t,x)|^2dx\right)dt
                             \medskip\\
                  \le & \;CA^2\vint_{I\cap\left(T_0+\frac{\kappa\eta\prime}{A^2},T_1-\frac{\kappa\eta\prime}{A^2}\right)}
                            \left(\int_{K(t)}|u(t,x)|^2dx\right)dt
                            \medskip \\
                        & +CA^2\|u\|_{L^\frac{2(N+2)}{N}((T_0,T_1)\times\R^N)}^2\left(\int_{I\cap\left[\left(T_0,
                            T_0+\frac{\kappa\eta\prime}{A^2}\right)\cup\left(T_1-
                            \frac{\kappa\eta\prime}{A^2},T_1\right)\right]}\left(\int_{K(t)}1\,
                            dx\right)dt\right)^\frac{2}{N+2}
                            \medskip \\
                   \le & \; CA^2|I|\vsup_{t\in I\cap\left(T_0+\frac{\kappa\eta\prime}{A^2},T_1
                            \frac{\kappa\eta\prime}{A^2}\right)}\int_{K(t)}|u(t,x)|^2dx
                            +CA^2\eta\prime^\frac{N}{N+2}\left(\frac{\kappa\eta\prime}{A^2}\right)^\frac{2}{N+2}
                            \bigg(\frac 1{A^2}\bigg)^{\frac{N}{N+2}} \medskip \\
                   \le & \;C\vsup_{t\in I\cap\left(T_0+\frac{\kappa\eta\prime}{A^2},
                            T_1-\frac{\kappa\eta\prime}{A^2}\right)}\int_{K(t)}|u(t,x)|^2dx
                            +C\kappa^\frac{2}{N+2}\eta\prime,
\end{align*}
where $C=C(\|u_0\|_{L^2},N,\eta).$ For such a $C,$ let $\kappa>0$ be small enough to have $C\kappa^\frac{2}{N+2}\le\frac{1}{2}.$ Then $\kappa=\kappa(\|u_0\|_{L^2},N,\eta)$ and
\begin{gather*}
\vsup_{t\in I\cap\left(T_0+\frac{\kappa\eta\prime}{A^2},T_1-\frac{\kappa\eta\prime}{A^2}\right)}\int_{K(t)}|u(t,x)|^2dx
\ge C\eta^\frac{2(N+2)}{N},
\end{gather*}
where $C=C(\|u_0\|_{L^2},N,\eta).$ So there exists $t_0\in I\cap\left(T_0+\frac{\kappa\eta\prime}{A^2},T_1-\frac{\kappa\eta\prime}{A^2}\right)$
such that
\begin{gather}
\label{dempropmass9}
\int_{K(t_0)}|u(t_0,x)|^2dx\ge C\eta^\frac{2(N+2)}{N},
\end{gather}
where $C=C(\|u_0\|_{L^2},N,\eta).$ Since $\ell(K(t_0))=\dfrac{1}{A},$ then $K(t_0)$ is contained in a ball of radius $\dfrac{\sqrt N}{A}.$ Furthermore, $T_0+\dfrac{\kappa\eta\prime}{A^2}<t_0<T_1-\dfrac{\kappa\eta\prime}{A^2},$ which yields
\begin{gather}
\label{dempropmass10}
\frac{1}{A}\le C\min\{(T_1-t_0)^\frac{1}{2}, (t_0-T_0)^\frac12\},
\end{gather}
where $C=C(\|u_0\|_{L^2},N,\eta).$ Using this and Step~2, it follows that $K(t_0)$ can be covered by a finite number (which depends only on $\|u_0\|_{L^2},$ $N$ and $\eta$) of balls of
radius $R=\min\left\{(T_1-t_0)^\frac{1}{2}, (t_0-T_0)^{\frac12}\right\}.$ Then, by (\ref{dempropmass9}), there is some $c\in\R^N$ such that
\begin{gather}
\int_{B(c,R)}|u(t_0,x)|^2dx\ge \eps(\|u_0\|_{L^2},N,\eta).
\end{gather}
This concludes the proof.
\medskip
\end{proof*}

\begin{vproof}{of Theorem~\ref{thmmass}.}
Let $\gamma,$ $u_0$ and $u$ be as in Theorem~\ref{thmmass}. Let $\eta_0=\eta_0(N,|\gamma|)>0$ be given by
Proposition~\ref{propmass}. We apply Proposition~\ref{propmass} with $\eta=\eta_0.$ Let $\eps=\eps(\|u_0\|_{L^2},N,|\gamma|)>0$ be given by Proposition~\ref{propmass}. Assume that $T_\M<\infty.$
Then $\|u\|_{L^\frac{2(N+2)}{N}((0,T_\M);L^\frac{2(N+2)}{N}(\R^N))}=\infty$
and so there exist
$$
0=T_1<T_2<\dots<T_n<T_{n+1}<\dots<T_\M
$$
such that
$$
\forall n\in\N,\; \|u\|_{L^\frac{2(N+2)}{N}((T_n,T_{n+1})\times\R^N)}=\eta_0.
$$
It follows from Proposition~\ref{propmass} that for each $n\in\N,$ there exist $c_n\in\R^N,$ $R_n>0$ and $t_n\in(T_n,T_{n+1})$ such that
$$
R_n\le \min\{(T_\M-t_n)^\frac{1}{2}, (T_\m+t_n)^\frac{1}{2} \}
\quad
\mbox{and} \quad
\|u(t_n)\|_{L^2(B(c_n,R_n))}\ge\eps,
$$
for every $n\in\N.$ The case $T_\m<\infty$ follows in the same way. Hence we have proved the result.
\medskip
\end{vproof}

\section{Further Results}
\label{further}

As a corollary of the previous results, we can generalize to higher dimensions the 2--dimensional results proved by Merle and Vega~\cite{MR1628235} and the results proved by Keraani in~\cite{MR2216444} dimensions 1 and 2. We state here the most interesting of them. We need first some notation.

\begin{defi}
Let $\gamma\in\R\setminus\{0\}.$ We define $\delta_0$ as the supremum of $\delta$ such that if
$$
\|u_0\|_{L^2}< \delta,
$$
then \eqref{nls} has a global solution $u\in C(\R;L^2(\R^N))\cap L^{\frac{2(N+2)}N}(\R;L^{\frac{2(N+2)}N}(\R^N)).$
\end{defi}

We can prove the following

\begin{thm}
\label{thmmasskeraa}
Let $\gamma\in\R\setminus\{0\},$ let $u_0\in L^2(\R^N)\setminus\{0\},$ such that
$\|u_0\|_{L^2(\R^N)}<\sqrt2\delta_0,$ and let
$$
u\in C((-T_\m,T_\M);L^2(\R^N))\cap L_\loc^\frac{2(N+2)}{N}((-T_\m,T_\M);L^\frac{2(N+2)}{N}(\R^N))
$$
be the maximal solution of $(\ref{nls})$ such that $u(0)=u_0.$ Assume that $T_\M<\infty,$ and let $\lambda(t)>0,$ such that $\lambda(t)\tends\infty$ as $t\tends T_\M.$ Then there exists $x(t)\in\R^N$ such that,
\begin{gather*}
\vliminf_{t\nearrow T_\M}\int_{B(x(t),\lambda(t)(T_\M-t)^\frac{1}{2})}|u(t,x)|^2dx\ge\delta_0^2.
\end{gather*}
If $T_\m<\infty$ and $\lambda(t)\tends\infty$ as $t\tends -T_\m$ then there exists $x(t)\in\R^N$ such that,
\begin{gather*}
\vliminf_{t\searrow -T_\m}\int_{B(x(t),\lambda(t)(T_\m+t)^\frac{1}{2})}|u(t,x)|^2dx\ge\delta_0^2.
\end{gather*}
\end{thm}

The main ingredient in the proof of that theorem is a profile decomposition of the solutions of the free Schr{\"o}dinger equation. This decomposition was shown in the case $N=2$ by Merle and Vega~\cite{MR1628235} (see also Theorem~1.4 in~\cite{MR2247881}) and by Carles and Keraani~\cite{MR2247881} when $N=1.$ We generalize it to higher dimensions thanks to the improved Strichartz estimate, Theorem~\ref{stric3}. To describe it we need a definition. We follow the notation of Carles and Keraani~\cite{MR2247881}.

\begin{defi}
If $\Gamma^j=(\rho_n^j,t_n^j,\xi_n^j,x_n^j)_{n\in\N},$ $j=1,2,\dots$ is a family of sequences in
$(0,\infty)\times\R\times \R^N\times \R^N,$ we say that it is an {\it orthogonal family} if for all $j\ne k,$
$$
\limsup_{n\rightarrow \infty}\left(\frac{\rho_n^j}{\rho_n^k}+\frac{\rho_n^k}{\rho_n^j}+\frac{|t_n^j-t_n^k|}{(\rho_n^j)^2}+\left|\frac{x_n^j-x_n^k}{\rho_n^j}+\frac{t_n^j\xi_n^j-t_n^k\xi_n^k}{\rho_n^j}\right|\right)=\infty.
$$
\end{defi}

Now, we can state the theorem about the linear profiles.

\begin{thm}
\label{thmlinearprofiles}
Let $(u_n)_{n\in\mathbb N}$ be a bounded sequence in $L^2(\R^N).$ Then, there exists a subsequence $($that we name $(u_n)$ for the sake of  simplicity$)$ that satisfies the following$:$ there exists a family $(\phi^j)_{j\in\N}$ of functions in $L^2(\R^N)$ and a family of pairwise orthogonal sequences $\Gamma^j=(\rho_n^j,t_n^j,\xi_n^j,x_n^j)_{n\in\N},$ $j=1,2,\dots$ such that
$$
\vT(t) u_n(x)=\sum_{j=1}^\ell H_n^j(\phi^j)(t,x)+w_n^\ell(t,x),
$$
where
$$
H_n^j(\phi)(t,x)=\vT(t)\left(e^{i(\cdot)\frac{\xi_n^j}2}\vT(-t_n^j)\frac1{(\rho_n^j)^{N/2}} \phi
\left(\frac{\cdot-x_n^j}{\rho_n^j}\right)\right)(x),
$$
with
$$
\limsup_{n\rightarrow\infty} \|w_n^\ell\|_{L^{\frac{2(N+2)}N}(\R\times\R^N)}\tends0\quad as \quad \ell\tends\infty.
$$
Moreover, for every $\ell\ge1,$
$$
\|u_n\|_{L^2(\R^N)}^2=\sum_{j=1}^\ell \| \phi^j\|_{L^2(\R^N)}^2+\|w_n^\ell(0)\|_{L^2(\R^N)}^2+{\mathrm o}(1),
$$
as $n\tends\infty.$
\end{thm}

A similar result has been proved for wave equations by Bahouri and G\'erard~\cite{MR1705001}. To prove Theorem~\ref{thmlinearprofiles} one can follow Carles and Keraani (proof of Theorem~1.4) in~\cite{MR2247881}. It is observed in that paper (Remark~3.5)  that the result follows from the refined Strichartz's estimate, our Theorem~\ref{stric3}, once we overcome a technical issue, due to the fact that the Strichartz exponent
$\frac{2(N+2)}N$ is an even natural number when $N\in\{1,2\}$ (which covers the cases that the previous authors considered) but not in higher dimensions (except $N=4).$ Thus, to complete the proof we only need the following orthogonality result.

\begin{lem}
\label{natural}
For any $M\ge1,$
$$
\|\sum_{j=1}^M H_n^j(\phi^j)\|_{L^{\frac{2(N+2)}N}(\R^{N+1})}^{\frac{2(N+2)}N}\le
\sum_{j=1}^M\|H_n^j(\phi^j)\|_{L^{\frac{2(N+2)}N}(\R^{N+1})}^{\frac{2(N+2)}N}+{\mathrm o}(1)
\quad as \quad n\tends\infty.
$$
\end{lem}

\begin{proof*}
The proof if based on a well-known orthogonality property (see G{\'e}rard~\cite{MR1632171} and (3.47) in Merle and Vega~\cite{MR1628235}): if we have two orthogonal families $\Gamma^1$ and $\Gamma^2,$ and two functions in $L^2(\R^N),$ $\phi^1$ and $\phi^2,$ then
\begin{equation}
\label{orthogonality}
\|H_n^1(\phi^1)H_n^2(\phi^2)\|_{L^{\frac{N+2}N}(\R^{N+1})}={\mathrm o}(1)\quad as \quad n\tends\infty.
\end{equation}
When $N=1$ or $N=2,$ $\frac{2(N+2)}N$ is a natural number, so we can decompose the $L^{\frac{2(N+2)}N}$ norm as a product and, using (\ref{orthogonality}), we obtain directly the lemma. In the
higher dimensional case, write
\begin{gather*}
 \begin{split}
     & \|\sum_{j=1}^MH_n^j(\phi^j)\|_{L^{\frac{2(N+2)}N}}^{\frac{2(N+2)}N}=\int|\sum_j
         H_n^j(\phi^j)|^2|\sum_j H_n^j(\phi^j)|^\frac{4}{N} \medskip \\
  = & \int\sum_{j}\sum_{k}|H_n^j(\phi^j)H_n^k(\phi^k)||\sum_\ell H_n^\ell(\phi^\ell)|^\frac{4}{N} \medskip \\
  = & \sum_j \int |H_n^j(\phi^j)|^2|\sum_\ell H_n^\ell(\phi^\ell)|^\frac{4}{N}+
         \sum_j\sum_{k\ne j}\int| H_n^j(\phi^j)H_n^k(\phi^k)||\sum_\ell H_n^\ell(\phi^\ell)|^\frac{4}{N} \medskip \\
  \stackrel{{\rm not}}{=} & A+B.
 \end{split}
\end{gather*}
We estimate $B$ using H{\"o}lder's inequality with exponents $\frac{N+2}N$ and $\frac{N+2}2,$
\begin{align*}
     & \int| H_n^j(\phi_j)H_n^k(\phi^k)||\sum_\ell H_n^\ell(\phi^j)|^\frac{4}{N}\\
\le & \|H_n^j(\phi^j)H_n^k(\phi^k)\|_{L^{\frac{N+2}N}(\R^{N+1})}\|\sum_{\ell=1}^M
          H_n^\ell(\phi^\ell)\|_{L^{\frac{2(N+2)}N}}^{\frac{4}N}.
\end{align*}
Then, we use the orthogonality \eqref{orthogonality} and obtain $B={\mathrm o}(1).$
\medskip \\
About $A,$ when $N\ge4$ then $\frac{4}{N}\le1$ and therefore,
\begin{gather*}
 \begin{split}
  A\le & \sum_j\sum_\ell\int|H_n^j(\phi^j)|^2|H_n^\ell(\phi^\ell)|^\frac{4}{N} \medskip \\
      =  & \sum_j \int|H_n^j(\phi^j)|^2|H_n^j(\phi^j)|^\frac{4}{N}+\sum_j\sum_{\ell\ne j}\int|H_n^j(\phi^j)|^2
              |H_n^\ell(\phi^\ell)|^\frac{4}{N}.
 \end{split}
\end{gather*}
The first term of the sum is
$$
\sum_j\|H_n^j(\phi^j)\|_{L^{\frac{2(N+2)}N}}^{\frac{2(N+2)}N}.
$$
The second one is
$$
\sum_j\sum_{\ell\ne j}\int|H_n^j(\phi^j)|^{2-\frac{4}{N}}|H_n^j(\phi^j)H_n^\ell(\phi^\ell)|^\frac{4}{N}.
$$
We apply H{\"o}lder's with exponents $\frac{N+2}{N-2}$ and $\frac{N+2}4$ and bound the last sum by
$$
\sum_j\sum_{j\ne \ell}\|H_n^j(\phi_n^j)|\|_{L^{\frac{2N+4}N}}^{2-\frac{4}{N}}
\|H_n^j(\phi^j)H_n^\ell(\phi^\ell)\|_{L^{\frac{N+2}N}}^{\frac{4}N}
$$
which is ${\mathrm o}(1)$ by \eqref{orthogonality}. This finishes the proof of the Lemma for $N\ge4.$
\medskip \\
When $N=3,$ then $\frac{4}{N}=\frac{4}{3}>1,$ which complicates a bit the argument. We write
\begin{gather*}
A=\sum_j \int |H_n^j(\phi_j)|^2|\sum_\ell H_n^\ell(\phi^\ell)||\sum_m H_n^m(\phi^m)|^\frac{1}{3}
\le\sum_\ell\sum_j\sum_m \int |H_n^j(\phi^j)|^2|H_n^\ell(\phi^\ell)||H_n^m(\phi^m)|^\frac{1}{3}.
\end{gather*}
Using a similar argument as in the previous case, we show that the above integrals are $\rm{o}(1)$ except in the case $j=\ell=m.$
This ends the proof of the lemma for $N=3.$
\medskip
\end{proof*}

\begin{vproof}{of Theorem~\ref{thmmasskeraa}.}
To prove Theorem~\ref{thmmasskeraa}, one can follow  the arguments given by Keraani in~\cite{MR2216444}. Again one has to deal with the fact that $\frac{4}{N}$ is not in general a natural number. Apart from Lemma~\ref{natural}, we just need an elementary inequality (see (1.10) in G{\'e}rard~\cite{MR1632171}) for the function $F(x)=|x|^\frac{4}{N}x:$
$$
|F(\sum_{j=1}^\ell U^j)-\sum_{j=1}^\ell F(U^j)|\le\sum_j\sum_{k\ne j}|U^j||U^k|^\frac{4}{N}.
$$
Then, the arguments given by Keraani generalize to higher dimensions without difficulty, and prove
Theorem~\ref{thmmasskeraa}.
\medskip
\end{vproof}

\begin{rmk}
\label{rmkkeraani}
As said in the beginning of this section, we generalize all the results of Keraani~\cite{MR2216444} to higher dimension $N.$ In
particular, we display two of them.
\begin{enumerate}
\item
There exists an initial data $u_0\in L^2(\R^N)$ with $\|u_0\|_{L^2}=\delta_0,$ for which the solution $u$ of
\eqref{nls}  blows-up in finite time $T_\M.$
\item
Let $u$ be a blow-up solution of  \eqref{nls} at finite time $T_\M$ with initial data $u_0,$ such that
$\|u_0\|_{L^2}<\sqrt2\,\delta_0.$ Let $(t_n)_{n\in\N}$ be any time sequence such that $t_n\xrightarrow[]{n\to\infty}T_\M.$ Then there exists a subsequence of $(t_n)_{n\in\N}$ (still denoted by
$(t_n)_{n\in\N}),$ which satisfies the following properties. There exist $\psi\in L^2(\R^N)$ with $\|\psi\|_{L^2}\ge\delta_0,$ and a sequence $(\rho_n,\xi_n, x_n)_{n\in\N}\in(0,\infty)\times\R^N\times\R^N$ such that
$$
\vlim_{n\to\infty}\frac{\rho_n}{\sqrt{T_\M-t_n}}\le A,
$$
for some $A\ge0,$ and
$$
\rho_n^\frac{N}{2}e^{ix\xi_n}u(t_n,\rho_n x+x_n)\weak\psi \; \mbox{ in } \; L^2_\w(\R^N),
$$
as $n\tends\infty.$
\end{enumerate}
\end{rmk}

\noindent
{\bf Acknowledgments.} The authors are grateful to Professors James Colliander and Patrick G{\'e}rard for their very useful comments. We also would like to thank Professor Sahbi Keraani for his willingness to show us his work~\cite{MR2216444}.

\baselineskip .0cm

\bibliographystyle{abbrv}
\bibliography{Paper6}
\addcontentsline{toc}{section}{References}

\end{document}